\documentclass[mathpazo]{cicp}

\usepackage{graphicx}
\usepackage[dvipsnames]{xcolor}

\usepackage{amssymb}
\usepackage{amsmath}
\usepackage{dsfont}
\usepackage{rotating}
\usepackage{multirow}
\usepackage{hyperref}

\DeclareMathAlphabet\mathbfcal{OMS}{cmsy}{b}{n}  

\newcommand{\quotew}[1]{``#1''}
\newcommand{\x}{\mathbf{x}}
\newcommand{\q}{\mathbf{q}}
\newcommand{\CFL}{\textnormal{CFL}}
\newcommand{\dt}{\Delta t}
\newcommand{\R}{\mathds{R}}
\newcommand{\f}{\mathbf{f}} 
\newcommand{\g}{\mathbf{g}}
\newcommand{\Q}{\mathbf{Q}} 
\newcommand{\F}{\mathbf{F}}
\newcommand{\M}{\mathbf{M}}

\newcommand{\K}{\mathbf{K}}
\newcommand{\n}{\mathbf{n}}
\newcommand{\T}{\mathbf{T}}
\newcommand{\xxi}{\boldsymbol{\xi}}
\newcommand{\V}{\mathbf{V}}

\renewcommand{\u}{\mathbf{u}}
\renewcommand{\v}{\mathbf{v}}

\newcommand{\rd}[1]{{\leavevmode\color{black} #1}}
\newcommand{\bl}{\textcolor{black}}
\newcommand{\gf}{\textcolor{black}}
\newcommand{\Rall}{\textcolor{black}}

\begin{document}
\title{Continuous finite element subgrid basis functions for Discontinuous Galerkin schemes on unstructured polygonal Voronoi meshes}


 \author[BDG]{Walter Boscheri\affil{1}\comma\corrauth, Michael Dumbser \affil{2}, Elena Gaburro\affil{3}}
 \address{\affilnum{1}\ Department of Mathematics and Computer Science, University of Ferrara, Via Machiavelli 30, 44121 Ferrara, Italy. \\
           \affilnum{2}\ Department of Civil, Environmental and Mechanical Engineering, University of Trento, Via Mesiano 77, 38123 Trento, Italy \\
       \affilnum{3}\ Inria, Univ. Bordeaux, CNRS, Bordeaux INP, IMB, UMR 5251, 200 Avenue de la Vieille Tour, 33405 Talence cedex, France}
 \emails{{\tt walter.boscheri@unife.it} (W.~Boscheri), {\tt michael.dumbser@unitn.it} (M.~Dumbser),
          {\tt elena.gaburro@inria.fr} (E.~Gaburro)}

\begin{abstract}
We propose a new high order accurate \textit{nodal} discontinuous Galerkin (DG) method for the solution of nonlinear hyperbolic systems of partial differential equations (PDE) on unstructured polygonal Voronoi meshes. Rather than using classical \textit{polynomials} of degree $N$ inside each element, in our new approach the discrete solution is represented by \textit{piecewise continuous polynomials} of degree $N$ \textit{within} each Voronoi element, using a \textit{continuous finite element} basis defined on a subgrid inside each polygon. 
We call the resulting subgrid basis an \textit{agglomerated finite element} (AFE) basis for the DG method on general polygons, since it is obtained by the agglomeration of the finite element basis functions associated with the subgrid triangles. 
The basis functions on each sub-triangle are defined, as usual, on a universal reference element, hence allowing to compute  \textit{universal} mass, flux and stiffness matrices for the subgrid triangles once and for all in a pre-processing stage for the reference element only.  Consequently, the construction of an efficient \textit{quadrature-free} algorithm is possible, despite the unstructured nature of the computational grid. 
High order of accuracy in time is achieved thanks to the ADER approach, making use of an element-local space-time Galerkin finite element predictor.

The novel schemes are carefully validated against a set of typical benchmark problems for the compressible Euler and Navier-Stokes equations. The numerical results have been checked with reference solutions available in literature and also systematically compared, in terms of computational efficiency and accuracy, with those obtained by the corresponding modal DG version of the scheme.
\end{abstract}


\maketitle


\section{Introduction} \label{sec.intro}
Nonlinear systems of hyperbolic conservation laws are used to model a wide range of phenomena in nature, covering different fields of applications in science and engineering. Due to the complex structure of the governing partial differential equations (PDE) and the non-linearity of the problem, analytical solutions are extremely rare and very difficult to be found. Therefore, a lot of research has been devoted to the design of numerical schemes that aim at solving nonlinear systems of evolutionary PDE. Suitable discretizations in both space and time have been investigated in the past decades, starting from the pioneering work on Godunov-type finite volume schemes~\cite{GodunovRS,vanLeerRS,ToroBook}. In this approach the numerical solution is stored under the form of piecewise \textit{constant} cell averages within each control volume of the computational mesh, thus requiring a spatial reconstruction procedure in order to obtain higher order schemes, and the time evolution is obtained either by using Runge-Kutta timestepping, or directly discretizing the integral form of the conservation law at the aid of a fully-discrete one-step method. Alternatively, discontinuous Galerkin (DG) finite element methods can be used for the spatial approximation of the numerical solution, that in this case is directly expressed through high order \textit{polynomials} within each control volume, allowing jumps of the discrete solution across element boundaries, leading thus to a natural high order piecewise polynomial data representation. Thus, DG schemes do not need any reconstruction procedure, unlike high order finite volume solvers. These methods were first applied to neutron transport equations~\cite{reed} and later extended to general nonlinear systems of hyperbolic conservation laws in one and multiple space dimensions~\cite{cbs0,cbs1,cbs2,cbs3,cbs4}. 

In the DG framework, the numerical solution is represented globally by piecewise polynomials of degree up to $N$ using a polynomial expansion in terms of a set of suitable basis functions inside each element, that can be either of nodal or modal type. The \textit{nodal approach} allows very efficient schemes to be formulated in terms of a piecewise polynomial approximation with a judicious choice of the nodes, like the Gauss-Lobatto nodes~\cite{GassnerDG_LES} or the Gauss-Legendre nodes~\cite{Exahype}. Nodal basis functions are typically used either on Cartesian meshes or on unstructured grids composed of simplex control volumes, namely triangles in 2D and tetrahedra in 3D, but there exist also nodal basis functions for more general unstructured meshes, see, e.g. \cite{GassnerPoly}. This restriction is mainly due to the preferable requirement of a reference element where the basis functions are uniquely defined and the physical control volume is mapped to. Hence, the nodal basis is universal and only the mapping between the physical and the reference coordinate system contains the element-dependent information. Consequently, mass, stiffness and flux matrices, that typically arise from a DG discretization of the governing PDE, can be conveniently computed only once and for all on the reference element (e.g. the unit square/cube or the reference triangle/tetrahedron in 2D/3D, respectively) and used throughout the entire computation for all cells and time steps. Furthermore, if the nodal points of the basis coincide with the nodes of a sufficiently accurate quadrature formula\footnote{For a nodal basis to be truly orthogonal, without any kind of mass lumping or under-integration, the quadrature formula associated with the nodes needs to be exact for at least polynomials of degree $2 N$ in order to integrate the product of two basis functions, and thus the element mass matrix, exactly over each control volume.}, such as the Gauss-Legendre points for example, then the degrees of freedom of the nodal basis automatically provide the value of the solution at the quadrature nodes and leads to an orthogonal basis by construction. This reduces the integral computation to a very efficient quadrature-free matrix-vector multiplication between the pre-computed universal integral matrix (evaluated on the reference element for the chosen nodal basis) and the vector of expansion coefficients of the numerical solution. 
On the other hand, also the \textit{modal approach} permits to deal with universal basis functions defined on the reference element~\cite{Dumbser2008,BalsaraDGMaxwell2019}. However, in this case, the expansion coefficients of the numerical solution are the modes of the polynomial basis and therefore they do not directly provide the values of the numerical solution at the quadrature-nodes. Nevertheless, quadrature-free schemes can still be constructed at the aid of modal basis functions, see~\cite{atkins,fambri2017space}. Moreover, the adoption of a hierarchical modal approach is useful in the design of limiting techniques for ensuring the stability of DG schemes, such as classical slope and moment limiters of the DG method, see e.g. \cite{cbs4,Biswas_94,Burbeau_2001,Kri07,Kuzmin2013}. For more recent subcell finite volume limiters of the DG method the reader is referred to \cite{DGLimiter1,Sonntag2,DGsubcell_Gassner,MunzDGFV,ALEDG}, while artificial viscosity limiters are discussed, for example, in \cite{PerssonAV,VegtAV,TavelliCNS,GassnerMunzAV}. Indeed, one advantage of the modal basis is that they are able to deal with general polygonal/polyhedral computational cells because the basis functions can be directly defined in the physical coordinate system. In particular, because of their simplicity and high versatility, rescaled Taylor monomials are often used as modal basis functions on general unstructured grids~\cite{VoronoiDivFree,tavelli2014high,CWENOBGK,xu2016new}. Although exhibiting the formal order of accuracy of the scheme and proving their effectiveness and robustness in the numerical solution of PDE systems, the computational \textit{cost} associated to these modal DG schemes is rather high compared to the nodal approach. Indeed, no integrals can be precomputed on any reference element and all quantities must be evaluated at each space-time quadrature point whenever required by the numerical scheme. For example, in~\cite{ArepoTN, GaburroUnified} a modal expansion based on rescaled Taylor monomials was used in the context of finite volume and DG schemes for handling moving Voronoi meshes with topology changes, and a considerable amount of computational time is spent for the numerical flux evaluation since it must be performed in each quadrature point again throughout the entire simulation. Likewise, in the context of kinetic  equations~\cite{CWENOBGK}, a profiling analysis shows that the numerical integration over general polygonal meshes is the computationally most expensive part of the entire algorithm. In particular, in the aforementioned references, integration over the control volumes is carried out by splitting the polygonal cells into \textit{simplexes}, like triangles in 2D, and by employing Gauss quadrature formulae over each sub-triangular subcell, which further increases the computational efforts.   	
The problem of numerical quadrature on general polygonal meshes, including non-convex polygons, becomes even more complex in the context of \textit{agglomerated} DG finite element schemes, see for instance \cite{BassiAgglomeration, BassiAgglomeration2}. 

In order to overcome the problem of quadrature and in order to introduce a suitable continuous nodal basis  on general unstructured polygonal control volumes it is important to mention the very recent and computationally highly efficient virtual element method (VEM) for the numerical discretization of PDEs on general polygonal/polyhedral elements, see   \cite{Veiga_VEMbasic2013,Veiga_VEMhighpolyhedral2017,Veiga_VEM_NS_polygonal2018,Artioli_VEMcurvilinear2020,Veiga_VEM_C1_2020}. The VEM is a generalization of classical conforming \textit{continuous} finite element methods to arbitrary polygonal/polyhedral elements that makes use of \textit{non-polynomial} basis functions that are only \textit{implicitly known}, such that the elementary stiffness matrix of the associated variational problem can be computed without actually knowing explicitly the basis functions, but just using the known degrees of freedom on the boundary of the elements. 

In this work we aim at introducing a new piecewise polynomial nodal basis based on continuous finite element subgrid basis functions into the DG framework for the solution of conservation laws on unstructured polygonal Voronoi meshes. The new basis functions will be defined on a subgrid composed of triangular simplex subcells, where a continuous finite element approach will be exploited to represent the discrete solution \textit{within} each polygonal DG element. On the other hand, a discontinuous finite element method will be devised on the main cells given by the Voronoi polygons. In other words, the new basis used in this paper is \textit{not polynomial} in each Voronoi element, but \textit{piecewise polynomial}. Our method being a DG scheme, the discrete solution is still allowed to \textit{jump} across element boundaries, while it is \textit{continuous inside} each polygonal element.   
Fully discrete one-step DG schemes will be formulated relying on the \rd{so-called} ADER approach \rd{(Arbitrary order DERivative Riemann problem)}, proposed originally by Toro et al.~\cite{mill,toro2,toro3,titarevtoro,toro2006derivative,BTVC16}, but using here the more recent variant of~\cite{Dumbser2008}, which, in contrast to classical Runge-Kutta (RK) DG methods, does not require an evaluation at each intermediate Runge-Kutta stage of the semi-discrete scheme. 
\rd{For further details on the ADER approach we refer to the recent paper \cite{Busto2020_Frontiers} where an historical overview on the topic is presented together with a detailed description, for both structured and unstructured meshes, of the variant employed in this work.}
Thus, the technique presented in this paper proposes new high order accurate quadrature free and fully-discrete one-step DG schemes on two-dimensional Voronoi meshes. 
Detailed quantitative comparisons in terms of accuracy and computational efficiency against the same algorithm that instead makes use of traditional modal basis functions~\cite{ArepoTN} shows that the novel method is not only more accurate on the same mesh, but also computationally more efficient.

The outline of this article is as follows. In Section~\ref{sec.pde}, the governing equations are introduced. Section~\ref{sec.numscheme} is devoted to present all the details regarding the proposed numerical method, while in Section~\ref{sec.test} we show numerical convergence rates up to fourth order of accuracy in space and time for a smooth problem as well as a wide set of benchmark test problems considering inviscid and viscous compressible flows.
\Rall{Also a reasoned numerical and theoretical justification of the appropriateness of using Voronoi tessellations instead of the underlying Delaunay triangulation is given at the beginning of Section~\ref{sec.test}.}
Finally, we give some concluding remarks and an outlook to future work in Section~\ref{sec.concl}.

\section{Governing equations} \label{sec.pde}
We consider a bounded domain $\Omega \in \R^d$, where $d=2$ indicates the number of space dimensions, which is defined by the space coordinates $\x=(x,y)$, while $t$ refers to the time coordinate. The governing equations are given by general nonlinear viscous conservation laws of the form
\begin{equation}
\label{eqn.PDE}
\frac{\partial \Q}{\partial t} + \nabla \cdot \F(\Q,\nabla \Q) = \mathbf{0}, \qquad \x \in \Omega \subset \mathds{R}^d, \quad t \in \mathds{R}_0^+, \quad \Q \in \Omega_{\Q} \subset \mathds{R}^\nu,     
\end{equation} 
with $\Q$ denoting the vector of conserved variables defined in the space of the admissible states $\Omega_{\Q}\subset \R^\nu$ and $\F(\Q,\nabla \Q)=\left( \f(\Q,\nabla \Q),\g(\Q,\nabla \Q) \right)$ representing the nonlinear flux tensor which depends on the state $\Q$ and its gradient~$\nabla \Q$. 
The compressible Navier-Stokes equations for a Newtonian fluid with heat conduction can be cast in the form~\eqref{eqn.PDE} with
\begin{equation}
\label{eqn.NSEterms}
\Q=\left( \begin{array}{c} \rho \\ \rho \v \\ \rho E \end{array} \right), \qquad \F(\Q,\nabla \Q) = \left( \f(\Q,\nabla \Q),\g(\Q,\nabla \Q) \right)=
\left( \begin{array}{c}  \rho \v \\ \rho \left(\v \otimes \v \right) + \boldsymbol{\sigma}(\Q,\nabla \Q) \\ \mathbf{v} \cdot (\rho E \mathbf{I} + \boldsymbol{\sigma}(\Q,\nabla \Q)) - \kappa \nabla T  \end{array} \right), 
\end{equation}
where the $d \times d$ identity matrix is denoted with $\mathbf{I}$. The fluid density and pressure are $\rho$ and $p$, respectively, while $\mathbf{v}=(u,v)$ denotes the velocity vector. $E$ represents the specific total energy
\rd{given as the sum between the specific internal energy $e$ (see \eqref{eqn.idealEOS}) and the specific kinetic energy, i.e. $ E = e + \frac{1}{2} ||\mathbf{v}||^2$.}

 The stress tensor $\boldsymbol{\sigma}(\Q,\nabla \Q)$ is given under Stokes hypothesis as
\begin{equation}
\boldsymbol{\sigma}(\Q,\nabla \Q) = \left(p + \frac{2}{3} \mu \nabla \cdot \mathbf{v} \right) \mathbf{I} - \mu \left( \nabla \mathbf{v} + \nabla \mathbf{v}^T \right),
\label{eqn.sigma}
\end{equation}
with $\mu$ representing the dynamic viscosity that is assumed to be constant. The temperature is indicated with $T$ and the heat conduction coefficient $\kappa$ is linked to the viscosity through the Prandtl number $Pr$ as follows:
\begin{equation}
\kappa = \frac{\mu \gamma c_v}{Pr}.
\label{eqn.kappa}
\end{equation}
The adiabatic index $\gamma=c_p/c_v$ denotes the ratio of the specific heats at constant pressure $c_p$ and at constant volume~$c_v$, the latter given by $c_v=R/(\gamma-1)$ with $R$ being the specific gas constant. A thermal equation of state $p=p(T,\rho)$ and a caloric equation of state $e=e(T,\rho)$ are required to close the system~\eqref{eqn.NSEterms}. The temperature is typically canceled from these two equations of state, yielding one single relation of the form $e=e(p,\rho)$, which will be adopted in this work. Specifically, an ideal gas is considered with the thermal and caloric equation of state given by
\begin{equation}
\frac{p}{\rho}=RT, \qquad e=c_v T.
\label{eqn.idealEOS_twoeqn}
\end{equation}
The temperature can be eliminated using both expressions in~\eqref{eqn.idealEOS_twoeqn}, thus leading to an equation of state (EOS) of the form $e(p,\rho)$, that is
\begin{equation}
e(p,\rho) = \frac{p}{(\gamma-1)\rho}.
\label{eqn.idealEOS}
\end{equation}

For the two-dimensional compressible Navier–Stokes equations the convective eigenvalues $\boldsymbol{\lambda}$ and the viscous eigenvalues $\boldsymbol{\lambda}_v$ are given by, see~\cite{ADERNSE},   
\begin{equation}
\boldsymbol{\lambda}=(|\v|+c, \, |\v|, \, |\v|, \, |\v|-c), \qquad \boldsymbol{\lambda}_v=\left( \frac{4}{3} \frac{\mu}{\rho}, \, \frac{\gamma \mu}{Pr\,\rho} \right),
\label{eqn.eigenval}
\end{equation}
with $c^2=\gamma R T$ being the sound speed.

\section{Numerical scheme} \label{sec.numscheme}
The governing equations~\eqref{eqn.PDE} are numerically solved relying on fully-discrete one-step schemes that make use of a discontinuous Galerkin (DG) discretization on unstructured Voronoi meshes and an ADER time stepping technique~\cite{toro3,titarevtoro,Dumbser2008,Torlo_Offner_DEC-ADER2021} able to provide high order of accuracy in both space and time in one single step. 
The details of the numerical method are provided hereafter.

\subsection{Discretization of the space-time computational domain} \label{ssec.domain}
Let us start by providing the details about the discretization of the spatial computational domain and the time interval.

\begin{figure}[!htbp]
	\begin{center}
		\begin{tabular}{cc}  
			\includegraphics[width=0.47\textwidth]{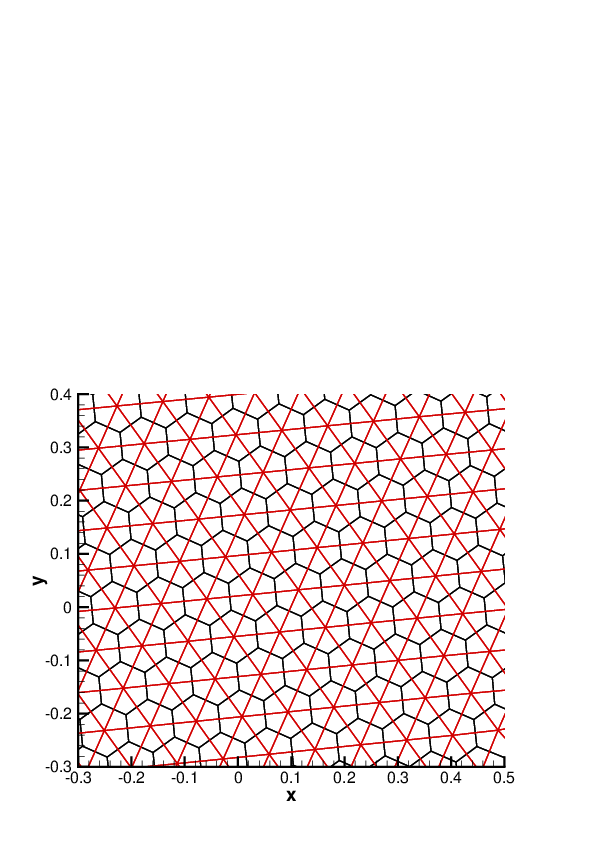} &
			\includegraphics[width=0.47\textwidth]{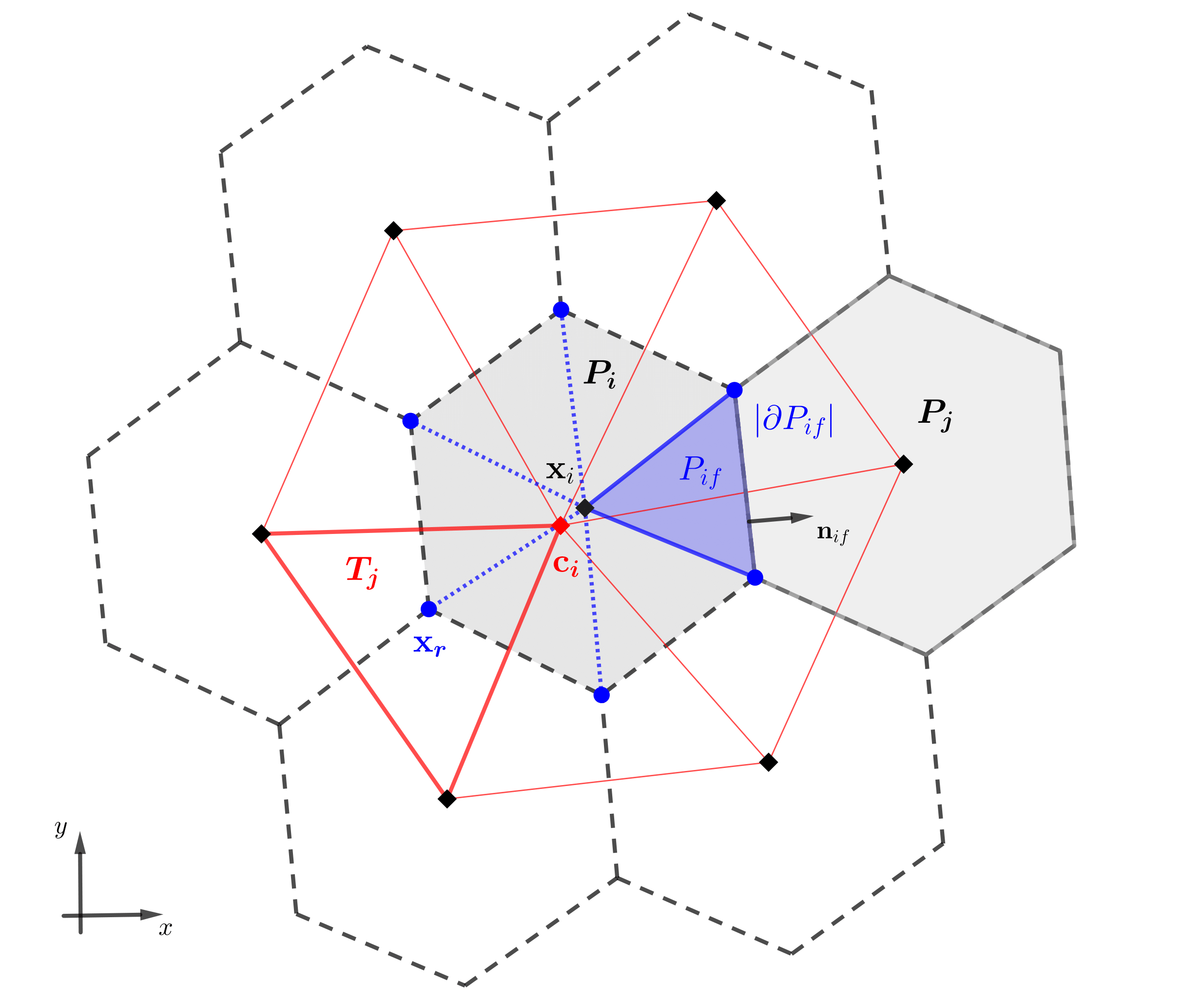} \\ 
		\end{tabular} 
		\caption{Example of Voronoi tessellation and mesh notation. Left: primary Delaunay triangulation (red lines) and dual Voronoi mesh (black lines). Right: Voronoi cell $P_i$ (filled in gray) with generator $\mathbf{c}_i$ (red dot), barycenter $\x_i$ (black dot), one of its vertex $\x_r \in \mathcal{R}_i$ (which coincides with the barycenter of a primary element $T_j$), face $\partial P_{if}$ and associated subcell $P_{if}$ (filled in blue). The outward pointing normal vector $\n_{if}$ and the Neumann neighbor $P_j$ are also depicted.}
		\label{fig.mesh}
	\end{center}
\end{figure}

\paragraph{Space discretization} The computational domain $\Omega$ is discretized with a centroid based Voronoi-type tessellation made of $N_E$ non-overlapping control volumes $P_i, \, i=1, \ldots, N_E$, which constitute the mesh configuration $\mathcal{D}_{\Omega}$:
\begin{equation}
\mathcal{D}_{\Omega} = \bigcup \limits_{i=1}^{N_E}{P_i}. 
\label{eqn.mesh}
\end{equation}
In order to construct the Voronoi mesh a total number of $N_E$ generator points with coordinates $\mathbf{c}_i, \, i=1, \ldots, N_E$ are needed: they are given by the vertexes of a \textit{primary} Delaunay triangulation that we have obtained by the mesh generator \textit{Gmsh}~\cite{GMSH}, which takes a characteristic target length $h$ as input parameter. The defining property of the Delaunay triangulation is that the circumcircle of each triangle $T_j$, $j=1, \ldots, N_T$ is not allowed to contain any of the other generator points in its interior. Each Voronoi element $P_i$ is then associated to a generator point $\mathbf{c}_i$ and assembled by connecting counterclockwise the barycenters of all the Delaunay triangles having this generator point as a vertex; see Figure~\ref{fig.mesh} for an example of a primary Delaunay mesh and the associated Voronoi tessellation. The usage of barycenters (instead of circumcenters) allows degenerate situations to be avoided, that might arise from the violation of the empty circumcircle property. This Voronoi tessellation is also referred to as the \textit{dual} mesh associated to the primary Delaunay triangulation.

In what follows we assume that the index $i$ refers to the element $P_i$, $f$ identifies a face (which is of dimension $d-1$) and $r$ denotes a vertex. The Neumann neighbor of element $P_i$ which shares with it the face $f$ is denoted with $P_j$, and the outward pointing unit normal vector is $\n_{if}$. At the aid of Figure~\ref{fig.mesh}, let us now introduce the following notation:
\begin{itemize}
	\item $\mathcal{R}_i$ is the set of $N_{R_i}$ vertexes $r$ of the Voronoi element $P_i$, whose coordinates $\x_{r}$ are given by the barycenters of the Delaunay triangles $T_j$, $j=1, \ldots, N_T$ of the primary triangulation;
	\item $\mathcal{R}_f$ is the set of vertexes $r$ of face $f$, which always accounts for two nodes in $d=2$;
	\item $\mathcal{F}_i$ is the set of $N_{R_i}$ faces $\partial P_{if}$ of the Voronoi element $P_i$; 
	\item $h_i$ denotes the characteristic mesh size of a Voronoi element $P_i$ which is defined as
	\begin{equation}
	h_i=\frac{2 \, |P_i|}{ \sum \limits_{f \in \mathcal{F}(P_i)} |\partial P_{if}|}, 
	\label{eqn.hi}
	\end{equation}
	where $|P_i|$ is the surface of the cell and $|\partial P_{if}|$ is the length of the face $\partial P_{if}$;
	\item $\x_i$ refers to the barycenter of a Voronoi element, that usually does not coincide with the generator point $\mathbf{c}_i$, and it is computed as
	\begin{equation}
	\x_{i} = \frac{1}{N_{R_i}} \sum \limits_{r \in \mathcal{R}_i} \x_r.
	\end{equation}
\end{itemize}
By connecting $\x_i$ with each vertex of $\mathcal{R}_i$, the Voronoi polygon $P_i$ is subdivided in $N_{R_i}$ subcells. Given a cell $P_i$ and one of its face $f \in \mathcal{F}_i$, the subcell $P_{if}$ is the triangle defined by three points, namely the barycenter $\x_i$ and the two vertexes of the set $\mathcal{R}_f$ associated to face $f$. Consequently, the index $if$ refers to the subcell $P_{if}$ and the Voronoi cell can also be seen as the union of its subcells, that is
\begin{equation}
P_i = \bigcup \limits_{f \in \mathcal{F}_i} P_{if}.
\label{eqn.subtri}
\end{equation}

\paragraph{Time discretization} The time interval $[0;t_f]$ is  discretized into time steps such that $t \in [t^n;t^{n+1}]$:
\begin{equation}
t = t^n + \tau \dt, \qquad \tau \in [0;1],
\label{eqn.time_map}
\end{equation}
with $t^n$ and $\dt:=(t^{n+1}-t^n)$ representing the current time and time step, respectively, while $\tau$ is the reference time coordinate which maps the time step $\dt$ to the reference interval $[0;1]$ according to~\eqref{eqn.time_map}. The size of the time step is determined by a CFL-type stability condition for explicit DG schemes which reads
\begin{equation}
\dt \leq \CFL \, \frac{\min \limits_{i \in N_E} h_i}{(2N+1) \max \limits_{i \in N_E} \left( |\lambda^{\max,i}| + 2 |\lambda_v^{\max,i}| \frac{2N+1}{h_i}\right)},
\label{eqn.timestep}
\end{equation}
\rd{where $N$ represents the degree of the chosen piecewise polynomial data representation, as described below in Section~\ref{ssec.numsol}.}

\subsection{Data representation on Voronoi meshes} \label{ssec.numsol}
A cell-centered discretization is adopted to represent the numerical solution for the conserved quantities $\Q$ in~\eqref{eqn.PDE} within each Voronoi polygon $P_i$ at the current time $t^n$, that is given in terms of piecewise polynomials of degree $N \geq 0$ denoted by $\mathbf{u}_h^n(\x,t^n)$ and defined in the space $\mathcal{U}_h$ as
\begin{equation}
\mathbf{u}_h^n(\x,t^n) = \sum \limits_{\rd{\ell=0}}^{\rd{\mathcal{N}_i-1}} \phi_\ell(\x) \, \hat{\mathbf{u}}^{n}_{\ell,i} 
:= \phi_\ell(\x) \, \hat{\mathbf{u}}^{n}_{\ell,i} , \qquad \x \in \rd{P_i},
\label{eqn.uh}
\end{equation}
where $\phi_\ell(\x)$ are the spatial basis functions used to span the space of (piecewise) polynomials $\mathcal{U}_h$ up to degree $N$ and $\mathcal{N}_i$ represents the total number of degrees of freedom for the cell $P_i$. For the sake of clarity, classical tensor index notation based on the Einstein summation convention is adopted, which implies summation over two equal indices. To complete the description of the expansion~\eqref{eqn.uh} the set of basis functions $\phi_\ell$ must be determined and this constitutes the main focus of the present work. 

\subsubsection{Modal basis functions}
The unstructured nature of the computational grid, where each Voronoi polygon might exhibit a different number of edges, makes the definition of a single reference element $P_e$ impossible for the general control volume $P_i$. 
\bl{Thus, a common solution consists in employing \textit{modal basis functions} i.e. functions directly defined in the \textit{physical space} with coordinates $\x$ for the definition of the numerical solution~\eqref{eqn.uh}.} 
A widespread technique is the usage of rescaled Taylor monomials of degree $N$ in the variables $\mathbf{x}=(x,y)$ given for the physical element $P_i$, expanded about its barycenter $\x_i$ and normalized by its characteristic length $h_i$, with $\ell=(\ell_1,\ell_2)$ being a multi-index:
\begin{equation} 
\label{eqn.modal_basis}
\phi_\ell(\x) |_{P_i} = \frac{(x - x_i)^{\ell_1}}{\ell_1! \, h_i^{\ell_1}} \, \frac{(y - y_i)^{\ell_2}}{\ell_2! \, h_i^{\ell_2}}, \qquad 
\ell = 0, \dots, \mathcal{N}-1, \quad \ 0 \leq \ell_1 + \ell_2 \leq N.
\end{equation} 
\rd{To ease the understanding of the above formula we remark that the space of polynomials of degree up to $N$ can be spanned by a set of modal basis functions containing $\mathcal{N}$ functions indexed by the symbol $\ell$. For each $\ell$ we uniquely define the values $\ell_1$ and $\ell_2$, so we can associate to any $\ell$ the multi-index $(\ell_1, \ell_2)$ and we have a total of $\mathcal{N}$ multi-index $\ell$.}

The unknown expansion coefficients $\hat{\mathbf{u}}^{n}_{\ell,i}$ in~\eqref{eqn.uh} are thus the normalized derivatives $h_i^{\ell_1+\ell_2}\frac{\partial^{\ell_1+\ell_2}}{\partial x^{\ell_1} \partial y^{\ell_2}} {\mathbf{Q}}(\mathbf{x}_i,t^n)$ appearing in the Taylor series expansion of $\Q$ about $\x_i$ and time $t^n$. The total number $\mathcal{N}_i$ of expansion coefficients (or degrees of freedom) $\hat{\mathbf{u}}^{n}_{\ell,i}$ depends only on the polynomial degree $N$ and is given by $\mathcal{N}_i = \mathcal{L}(N,d)$, with 
\begin{equation}
\mathcal{L}(N,d) = \frac{1}{d!} \prod \limits_{m=1}^{d} (N+m),
\label{eqn.nDOF}
\end{equation}
with $d=2$ because only two-dimensional domains are considered in this paper. Notice that $\mathcal{N}_i$ is the same for all elements $P_i, \, i=1,\ldots,N_E$ and corresponds to the \textit{minimum} number of nodes necessary to define a polynomial of degree $N$ in $d$ space dimensions. 

\bl{This approach based on modal basis functions has been already adopted in the context of free surface flows~\cite{Voronoi,VoronoiDivFree,ADERFSE}, and more recently for kinetic equations~\cite{CWENOBGK,FVBoltz,DGBoltz}.
Also, moving meshes with topology changes have been considered in~\cite{ArepoTN} for applications in fluid dynamics.
The cited references demonstrate that both high order of accuracy and robustness can be achieved with this choice.}

However, the major drawback of modal basis functions is related to the associated computational efficiency, which is quite compromised as the order of accuracy of the scheme increases. Specifically, integration over each Voronoi polygon $P_i$ is numerically carried out by summing up the contribution of each subcell $P_{if} \in \mathcal{F}_i$ which in turn are computed via Gauss quadrature formulae of suitable order of accuracy~\cite{stroud} and thus with an increasing number of quadrature points.  Therefore, the evaluation of all the basis functions~\eqref{eqn.modal_basis} at each quadrature point for each control volume is extremely expensive from the computational viewpoint. In the case of DG schemes, where volume and boundary integrals must be evaluated at each time step, this disadvantage becomes even more evident~\cite{DGBoltz}. To overcome this problem, the vector of the basis functions evaluated at each Gauss point of all Voronoi elements could be stored once and for all in the pre-processing stage. Although this remedy undoubtedly reduces the computational efforts, this strategy might become prohibitive in terms of memory requirements, as the polynomial degree $N$ of the basis functions gets higher or the computational mesh undergoes refinements.

\subsubsection{Agglomerated subgrid Finite Element (AFE) basis functions} \label{ssec.AFE}
The aforementioned efficiency problem can be solved by the usage of a set of \textit{nodal basis functions} to represent the numerical solution given by~\eqref{eqn.uh}. However, since no canonical reference elements, like triangles or quadrilaterals, fit the general polygons of an unstructured Voronoi mesh, the design of a suitable nodal basis with the \textit{minimum} number of nodes $\mathcal{L}(N,d)$ that are necessary to define a polynomial of degree $N$ in $d$ space dimensions is not straightforward. 

\begin{figure}[!htbp]
	\begin{center}
		\begin{tabular}{c}  
			\includegraphics[width=0.9\textwidth]{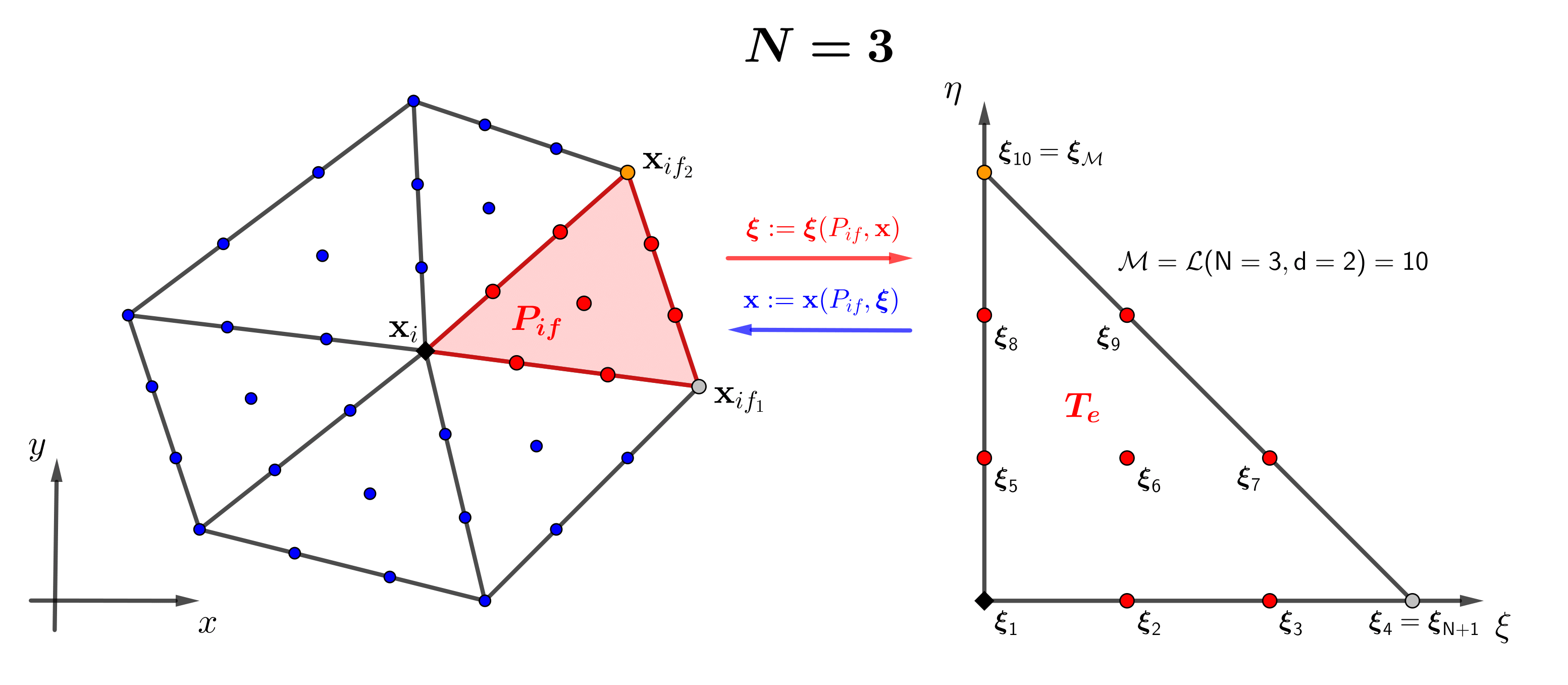} \\ 
		\end{tabular} 
		\caption{Agglomerated finite element (AFE) basis functions on the Voronoi polygon $P_i$ for degree $N=3$. Left: degrees of freedom on the sub-triangulation in the physical coordinate system $\x$ and subcell $P_{if}$ highlighted in red. Right: reference element $T_e$ in the coordinate system $\xxi$ to which the subcell $P_{if}$ is mapped  with corresponding local coordinates of the expansion coefficients. The barycenter $\x_i$ is mapped to the reference node $\xxi_1$, while $\x_{if_1}$ and $\x_{if_2}$ are mapped to the reference nodes $\xxi_{N+1}$ and $\xxi_{\mathcal{M}}$, respectively.}
		\label{fig.Te}
	\end{center}
\end{figure}

Hence, the main idea of this paper is to use the continuous union of piecewise polynomials of degree $N$ defined on each of the $N_{R_i}$ sub-triangles $P_{if}$ which constitute the Voronoi cell $P_i$ according to~\eqref{eqn.subtri}. The new basis functions are built following the standard nodal approach of classical \textit{conforming continuous finite elements} with order $N$ inside each subcell $P_{if}$, yielding a basis with a \textit{local} (i.e. within each sub-triangle) number of $\mathcal{M}=\mathcal{L}(N,d)$ degrees of freedom (DoF) given by~\eqref{eqn.nDOF}, but where the common DoF along the common edges of the subgrid triangles are \textit{shared}. 
In other words, the basis functions \textit{within} each DG element are given by the basis functions associated with classical continuous finite elements on a triangular subgrid defined within each polygon and as such is \textit{continuous inside} each polygon. But, as usual in the DG context, the discrete solution is \textit{allowed to jump across the boundary} $\partial P_i$ of $P_i$. 
Each subcell triangle can be easily mapped to the reference triangular element $T_e$ in the reference coordinate system $\xxi=(\xi,\eta)$ defined as $T_e = \{ (\xi,\eta) \in \R^2 : 0 \leq \xi \leq 1, \, 0 \leq \eta \leq 1-\xi  \}$. The spatial mapping of the triangular subcell $P_{if}$ from reference coordinates $\xxi$ to physical coordinates $\x$ relies on the simple linear coordinate transformation 
\begin{equation} 
\mathbf{x} := \x(P_{if},\xxi) = (1-\xi-\eta) \, \mathbf{x}_{i} + \xi \, \mathbf{x}_{if_1} + \eta \, \mathbf{x}_{if_2}, \qquad (\mathbf{x}_{if_1},\mathbf{x}_{if_2}) \in \mathcal{R}_f,
\label{eqn.xietaTransf} 
\end{equation} 
with $(\mathbf{x}_{if_1},\mathbf{x}_{if_2})$ denoting the two counterclockwise oriented vertexes of $P_i$ that define the face $f$ of the sub-triangle $P_{if}$, see Figure~\ref{fig.Te}. Notice that the barycenter $\x_i$ is always mapped to the first reference node at position $\xxi_1=(0,0)$. The inverse of the mapping is denoted by $\xxi:=\xxi(P_{if},\x)$, which directly results from inverting the linear relation~\eqref{eqn.xietaTransf}. Furthermore, the surface of the reference triangle is $|T_e|=1/2$, and the Jacobian matrix of the transformation for the subcell $P_{if}$ writes
\begin{equation}
\mathbf{J}_{if}= \left. \frac{\partial \x}{\partial \xxi} \right|_{P_{if}},
\end{equation}
with $|J_{if}|=2 \, |P_{if}|$ representing the determinant. For a visual interpretation, we refer to Figure~\ref{fig.Te} where we show the transformation~\eqref{eqn.xietaTransf} for $N=3$ with $\mathcal{M}=\mathcal{L}$$(N=3,$ $d=2)=10$ degrees of freedom in the reference element. 	

The coordinates of the nodal degrees of freedom are given in the reference element by
\begin{equation}
\xxi_k = \left( \xi_{k_1}, \eta_{k_2} \right) = \left(\frac{k_1}{N},\frac{k_2}{N}\right), \qquad 0 \leq k_1 \leq N, \qquad 0 \leq k_2 \leq (N-k_1),
\label{eqn.xinodes}
\end{equation}
with $k=(k_1,k_2)$ being a multi-index. The basis functions $\phi_\ell$ are then defined by the \textit{Lagrange interpolation polynomials} passing through the nodes $\xxi_k$, thus they satisfy the interpolation property 
\begin{equation}
\phi_\ell(\xxi_k) = \delta_{\ell k},
\label{eqn.delta_property}
\end{equation}
where $\delta$ indicates the Kronecker delta. This procedure completely defines the construction of the nodal basis functions over the subcell $P_{if}$.

To obtain the set of basis functions on the entire control volume $P_i$, the basis functions which are built over the subcell $P_{if}$ are trivially extended by zero to the rest of the element,
in order to maintain the continuity inside the Voronoi cell $P_i$. This yields the \textit{agglomerated finite element} (AFE) basis for each Voronoi element with the \textit{total} number of degrees of freedom $\mathcal{N}_i$ for $P_i$ that is given by
\begin{equation}
\mathcal{N}_i = N_{R_i} \left( \mathcal{L}(N,d) - N - 1 \right) + 1,
\label{eqn.spacedof} 
\end{equation}
since the degrees of freedom lying along the internal faces of the subcells are shared between the two neighboring sub-triangles and thus must be taken into account only once. The numerical solution~\eqref{eqn.uh} is then given by a total number of $\mathcal{N}_i$ degrees of freedom $\hat{\mathbf{u}}^{n}_{\ell,i}$. 
\bl{The degrees of freedom $\hat{\mathbf{u}}^{n}_{\ell,i}$ are given by the union, \textit{without repetition}, of the $\hat{\mathbf{u}}^{n}_{\ell,if}$ which are the $\mathcal{M}=\mathcal{L}(N,d)$ expansion coefficients related to the subcells $P_{if}$}, as shown in Figure~\ref{fig.Te} for the red sub-triangle.
Notice that in this case the total number of expansion coefficients $\mathcal{N}_i$ depends on both the polynomial degree $N$ and the number of faces of the Voronoi cell $P_i$, hence it is \textit{element-dependent}.

Compared to the modal basis~\eqref{eqn.modal_basis}, which accounts only for the minimum necessary number of degrees of freedom $\mathcal{L}(N,d)$ needed to achieve the formal order $N$ in $d$ space dimensions, the choice of the agglomerated finite element basis may not seem to be convenient  
at a first glance, since it involves a significantly higher number of expansion coefficients $\hat{\mathbf{u}}^{n}_{\ell,i}$ to represent the discrete solution.	
Nevertheless, these additional degrees of freedom carry also additional information, which provides more \textit{effective resolution} in the final numerical scheme, but without affecting negatively the computational cost, which instead is even reduced. 
Indeed, the main advantage of this new AFE subgrid basis for DG schemes is that the reference element can be used to map each subcell, hence allowing mass, flux and stiffness matrices to be computed only once in the reference system with Gauss quadrature rules of sufficient order to be exact and which are subsequently used to \textit{assemble} the final basis functions over the entire cell. 
In particular, by using the reference time $\tau$ of~\eqref{eqn.time_map} and the coordinate transformation~\eqref{eqn.xietaTransf}, the governing PDE~\eqref{eqn.PDE} can be written in the reference space-time coordinates $\tilde{\xxi}=(\xxi,\tau)$ as
\begin{equation}
\frac{\partial \Q}{\partial \tau} + \nabla_{\xxi} \cdot \F^*(\Q,\nabla \Q) = \mathbf{0}, \qquad \nabla_{\xxi} = \left( \frac{\partial}{\partial \xi}, \, \frac{\partial}{\partial \rd{\eta}} \right), \qquad \F^*(\Q,\nabla \Q) = (\f^*,\g^*),
\label{eqn.PDExi}
\end{equation}
with the transformed fluxes explicitly given by
\begin{equation}
\f^* = \dt \left( \f \frac{\partial \xi}{\partial x} + \g \frac{\partial \xi}{\partial y} \right), \qquad \g^{*} = \dt \left( \f \frac{\partial \eta}{\partial x} + \g \frac{\partial \eta}{\partial y} \right).
\label{eqn.fstar}
\end{equation} 	
This approach will ultimately lead to the design of very efficient \textit{quadrature-free} one-step DG schemes that are described in the next sections.

\subsection{Local space-time predictor} \label{ssec.predictor}
As introduced before, we will work in the context of fully discrete one-step ADER DG schemes which are generally constructed from two basic building blocks.
The first one, described in details in this section, consists in the computation of a \textit{local space-time predictor} solution $\q_h(\x,t)$ starting from the known numerical solution $\u_h(\x,t)$ at the current time $t^n$ for each considered cell $P_i$. This procedure respects the causality principle accounting only for information coming from the past in each cell, does not need any data communications between neighbor cells, and provides a high order space-time approximation of the discrete solution valid \textit{locally} inside $P_i \times [t^n;t^{n+1}]$, see \cite{eno,Dumbser2008}. 
Then, these predictor solutions will be used in the second part of the algorithm, the so-called corrector step (carefully described in the next Section~\ref{ssec.corrector}), where the actual solution $\hat{\mathbf{u}}^{n+1}_{\ell,i}$ at the new time $t^{n+1}$ will be finally recovered through an explicit fully discrete scheme obtained by integrating the governing PDE over $P_i \times [t^n;t^{n+1}]$. This explicit update formula will also account for the interaction between the predictors of neighbor cells thanks to the computation of Riemann solvers at the cell interfaces.

We refer the reader to \cite{ArepoTN, GaburroUnified} for the precise description of the predictor and corrector steps in the case of the modal basis functions~\eqref{eqn.modal_basis}, and we concentrate instead on the new agglomerated finite element basis introduced in Section ~\ref{ssec.AFE}. 
Thus, here, the predictor solution is represented within each \textit{subcell} $P_{if}$ by a space-time expansion of the form
\begin{equation}
\q_h(\x,t) = \sum \limits_{\ell=1}^{\mathcal{T}} \theta_\ell(\tilde{\xxi}) \hat \q_{\ell,if}:= \theta_\ell \hat \q_{\ell,if}, \qquad \x \in P_{if}, \qquad \mathcal{T}=\mathcal{L}(N,d+1),
\label{eqn.qh}
\end{equation}
with the unknown expansion coefficients $\hat \q_{\ell,if}$. These are a \textit{subset} of the expansion coefficients $\hat \q_{\ell,i}$ of the entire Voronoi cell, namely they belong only to the subcell $P_{if}$. The space-time basis functions $\theta_\ell=\theta_\ell(\tilde{\xxi})$ are defined by means of a tensor product between the agglomerated finite element basis in space and the one-dimensional Lagrange interpolation basis functions along the time coordinate passing through the Gauss-Legendre points \rd{(see Appendix A in \cite{BoscheriWAO} for explicit formulae of the basis functions on simplex control volumes)}. This yields the following set of space-time nodes for the subcell $P_{if}$ with the multi-index $k=(k_1,k_2,k_3)$:
\begin{equation}
\tilde{\xxi_k} = \left( \xi_{k_1}, \eta_{k_2}, \tau_{k_3} \right) = \left(\frac{k_1}{N},\frac{k_2}{N},\tau_{k_3}\right), \qquad 0 \leq k_1 \leq N, \quad 0 \leq k_2 \leq (N-k_1), 
\label{eqn.tildexinodes}
\end{equation}
with $\tau_{k_3}$ being the $k_3$-th root of the Legendre polynomial of degree $N$ rescaled to the unit interval $[0;1]$, and $\tilde{\xxi}=(\xxi,\tau)$ representing the space-time coordinates where the reference element $\tilde{T}_e = T_e \times [0;1]$ is defined. The total number of space-time degrees of freedom $\mathcal{N}_i^{st}$ for the element $P_i$ is
\begin{equation}
\mathcal{N}_i^{st} = \mathcal{N}_i \cdot  (N+1),
\end{equation}
with $\mathcal{N}_i$ given by~\eqref{eqn.spacedof}, thus the vector of unknown coefficients $\hat \q_{\ell,i}$ for the entire cell $P_i$ counts $\mathcal{N}_i^{st}$ degrees of freedom, which are collected from the expansion coefficients $\hat \q_{\ell,if}$ defined on the subcells $P_{if}, \, f \in \mathcal{F}_i$. 

The ADER approach is based on the solution of the generalized Riemann problem, which requires the knowledge of the time derivatives to evolve the solution in time. System~\eqref{eqn.PDE} is therefore solved \quotew{in the small} neglecting interactions between neighbor cells by relying on an element-local \textit{weak formulation} of the governing PDE~\eqref{eqn.PDExi} in the reference space-time coordinates. 
\rd{Indeed, the PDE in \eqref{eqn.PDExi} is integrated in space and time over the reference space--time element $\tilde{T}_e$ against a set of test functions $\theta_k(\tilde{\xxi})$ of the same form of the basis functions $\theta_l(\tilde{\xxi})$
\begin{equation} \nonumber 
	\int \limits_0^1 \int \limits_{T_e} \theta_k \frac{\partial \Q}{\partial \tau} \, d\tilde{\xxi} = - \int \limits_0^1 \int \limits_{T_e}  \theta_k \, \nabla_{\xxi} \cdot \F^*(\Q,\nabla \Q) \, d\tilde{\xxi}, 	\qquad \forall k \in [0, \mathcal{T}].
\end{equation}
Then, the unknown discrete solution $\q_h$ approximated by \eqref{eqn.qh} and the discrete flux $\F_h^*=(\f_h^*,\g_h^*)$ defined as 
\begin{equation} \nonumber
	\f_h^* = \sum \limits_{\ell=1}^{\mathcal{T}} \theta_\ell(\tilde{\xxi}) \hat \f_{\ell,if}^*:= \theta_\ell \hat \f_{\ell,if}^*, \qquad 
	\g_h^* = \sum \limits_{\ell=1}^{\mathcal{T}} \theta_\ell(\tilde{\xxi}) \hat \g_{\ell,if}^*:= \theta_\ell \hat \g_{\ell,if}^*, 
\end{equation}
are inserted in the weak form of the PDE replacing $\Q$ and $\F^*$, hence obtaining
\begin{equation} \nonumber 
	\int \limits_0^1 \int \limits_{T_e} \theta_k \frac{\partial \q_h}{\partial \tau} \,  d\tilde{\xxi} = - \int \limits_0^1 \int \limits_{T_e}  \theta_k \left ( \frac{\partial \f_h^*}{\partial \xi} + \frac{\partial \g_h^*}{\partial \eta} \right ) \, d\tilde{\xxi}, 	\qquad \forall k \in [0, \mathcal{T}],
\end{equation}
which can be expanded as
\begin{equation} \nonumber 
	\int \limits_0^1 \int \limits_{T_e} \theta_k \frac{\partial \theta_l}{\partial \tau} \, \hat \q_{\ell,if} \, d\tilde{\xxi} = - \int \limits_0^1 \int \limits_{T_e}  \theta_k \left ( \frac{\partial \theta_\ell}{\partial \xi} \, \hat \f_{\ell,if}^* +   \frac{\partial \theta_\ell}{\partial \eta} \, \hat \g_{\ell,if}^* \right ) \, d\tilde{\xxi}, 	\qquad \forall k \in [0, \mathcal{T}].
\end{equation}
Finally, remarking that the above expression holds true for all the sub-triangles $f \in \mathcal{F}_i$ of $P_i$, multiplication by $|J_{if}|$ and sum over $f$ lead to
\begin{equation} \nonumber 
	\sum \limits_{f\in \mathcal{F}_i} \, |J_{if}| \int \limits_0^1 \int \limits_{T_e} \theta_k \frac{\partial \theta_l}{\partial \tau} \, \hat \q_{\ell,if} \, d\tilde{\xxi} = - 
	\sum \limits_{f\in \mathcal{F}_i} \, |J_{if}| 
	\int \limits_0^1 \int \limits_{T_e}  \theta_k \left ( \frac{\partial \theta_\ell}{\partial \xi} \, \hat \f_{\ell,if}^* +   \frac{\partial \theta_\ell}{\partial \eta} \, \hat \g_{\ell,if}^* \right ) \, d\tilde{\xxi}, 	\qquad \forall k \in [0, \mathcal{T}].
\end{equation}
In order to uniquely determine each degree of freedom $\hat \q_{\ell,i}$ describing the predictor $\q_h$ over the Voronoi element $P_i$, on the left hand side we \textit{identify} the degrees of freedom $\hat q_{\ell, if}$ which refer to the same internal edges of the sub-triangles of $P_i$ without repetition and we get
\begin{equation} \nonumber 
	|P_i| \sum \limits_{f\in \mathcal{F}_i} \,  \int \limits_0^1 \int \limits_{T_e} \theta_k \frac{\partial \theta_l}{\partial \tau} \, \hat \q_{\ell,i} \, d\tilde{\xxi} = - 
	\sum \limits_{f\in \mathcal{F}_i} \, |J_{if}| 
	\int \limits_0^1 \int \limits_{T_e}  \theta_k \left ( \frac{\partial \theta_\ell}{\partial \xi} \, \hat \f_{\ell,if}^* +   \frac{\partial \theta_\ell}{\partial \eta} \, \hat \g_{\ell,if}^* \right ) \, d\tilde{\xxi}, 	\qquad \forall k \in [0, \mathcal{T}],
\end{equation}
which} compactly writes as
\begin{equation}
\K_{t,i} \hat \q_{\ell,i} = -  \, \sum \limits_{f \in \mathcal{F}_i} \rd{|J_{if}|} \left( \K_{\xi} \hat \f^{*}_{if} + \K_{\eta} \hat \g^{*}_{if} \right),
\label{eqn.PDEweak}
\end{equation}
where the following matrix definitions have been introduced:
\begin{equation}
\K_{\xi} = \int \limits_0^1 \int \limits_{T_e} \theta_k \frac{\partial \theta_\ell}{\partial \xi} \, d\tilde{\xxi}, \qquad \K_{\eta} = \int \limits_0^1 \int \limits_{T_e} \theta_k \frac{\partial \theta_\ell}{\partial \eta} \, d\tilde{\xxi}, \qquad \K_{t,i} = |P_i| \, \sum \limits_{f\in \mathcal{F}_i} \int \limits_0^1 \int \limits_{T_e} \theta_k \frac{\partial \theta_\ell}{\partial \tau} \, d\tilde{\xxi}.
\label{eqn.mat_def}
\end{equation}
Furthermore, the fluxes $\f^{*}_{if}$ and $\g^{*}_{if}$ are approximated using the same space-time basis functions $\theta_\ell(\tilde{\xxi})$ adopted for the numerical solution. Because of the interpolation property~\eqref{eqn.delta_property}, the expansion coefficients for the fluxes in~\eqref{eqn.PDEweak} can be directly evaluated in a \textit{pointwise} manner from $\hat \q_{\ell,if}$ using the definition~\eqref{eqn.fstar} as
\begin{equation}
\hat \f_{if}^* = \f^*(\hat \q_{\ell,if},\rd{\nabla \theta_\ell} \, \hat \q_{\ell,if}), \qquad \hat \g^*_{if} = \g^*(\hat \q_{\ell,if},\rd{\nabla \theta_\ell} \, \hat \q_{\ell,if}).
\end{equation}
Let us notice that the \textit{universal} stiffness matrices $\K_{\xi}$ and $\K_{\eta}$ are defined in the reference space-time element $\tilde{T}_e$, thus they do not depend neither on space nor on time, because the space-time dependency is taken into account by the space-time Jacobian $\dt |J_{if}|$ of the transformation from physical to reference coordinates. On the other hand, the time stiffness matrix $\K_{t,i}$ must handle simultaneously all degrees of freedom $\mathcal{N}_i^{st}$ for the predictor solution $\hat \q_{\ell,i}$, therefore it is \textit{element-dependent} and is of dimension $\mathcal{N}_i^{st} \times \mathcal{N}_i^{st}$. 

To introduce a proper upwinding approximation in the time direction, the term on the left hand side of~\eqref{eqn.PDEweak} is integrated by parts in time, which allows to take into account the initial condition of the local Cauchy problem in a weak form as 
follows: 
\begin{equation}
\K_{1,i} \hat \q_{\ell,i} = \sum \limits_{f \in \mathcal{F}_i} \rd{|J_{if}|} \, \F_0 \, \hat{\mathbf{u}}^{n}_{\ell,if} -  \sum \limits_{f \in \mathcal{F}_i} \rd{|J_{if}|} \, \left( \K_{\xi} \hat \f^{*}_{P_{if}} + \K_{\eta} \hat \g^{*}_{P_{if}} \right), 
\label{eqn.PDEweak2}
\end{equation}
with the matrices 
\begin{equation}
\K_{1,i} = |P_i| \,\sum \limits_{f\in \mathcal{F}_i}  \left[ \int \limits_{T_e} \theta_k(\xxi,1) \theta_\ell(\xxi,1) \, d\xxi - \int \limits_0^1 \int \limits_{T_e}  \frac{\partial \theta_k}{\partial \tau} \theta_\ell \, d\tilde{\xxi} \right], \qquad \F_0 = \int \limits_{T_e} \theta_k(\xxi,0) \phi_\ell(\xxi) \, d\xxi.
\label{eqn.mat_def2}
\end{equation}
Here, $\F_0$ is a universal matrix defined on the reference element $T_e$, while $\K_{1,i}$ is an element-wise matrix like the time stiffness matrix $\K_{t,i}$ in~\eqref{eqn.mat_def}. The above expression~\eqref{eqn.PDEweak2} constitutes an element-local nonlinear algebraic equation system for the unknown space-time expansion coefficients $\hat \q_{\ell,i}$ that is conveniently solved using the following iterative scheme:
\begin{equation}
\hat \q_{\ell,i}^{m+1} = \K_{1,i}^{-1} \left(  \sum \limits_{f \in \mathcal{F}_i} \rd{|J_{if}|} \, \F_0 \, \hat{\mathbf{u}}^{n}_{\ell,if} -  \sum \limits_{f \in \mathcal{F}_i} \rd{|J_{if}|} \,\left( \K_{\xi} \hat \f^{*,m}_{P_{if}} + \K_{\eta} \hat \g^{*,m}_{P_{if}} \right) \right),
\label{eqn.Picard}
\end{equation}
where $m$ denotes the iteration number. The iteration stops when the residual of~\eqref{eqn.Picard} is less than a prescribed tolerance, typically set to $10^{-12}$. \rd{We also remark that the convergence of the iterative solver has been demonstrated in~\cite{jackson2017eigenvalues} for linear PDEs and in~\cite{Busto2020_Frontiers} for nonlinear PDEs.}

To improve the computational efficiency, the inverse matrix $\K_{1,i}^{-1}$ can be computed and stored for each physical element $P_i$ in the pre-processing stage. This is the only relevant memory requirement of the predictor step because all other matrices in the nonlinear equation~\eqref{eqn.PDEweak2} are computed on the reference element where the novel agglomerated finite element basis functions are uniquely defined for each subcell $P_{if}$.

\begin{figure}[!htbp]
	\begin{center}
		\begin{tabular}{c}  
			\includegraphics[width=0.9\textwidth]{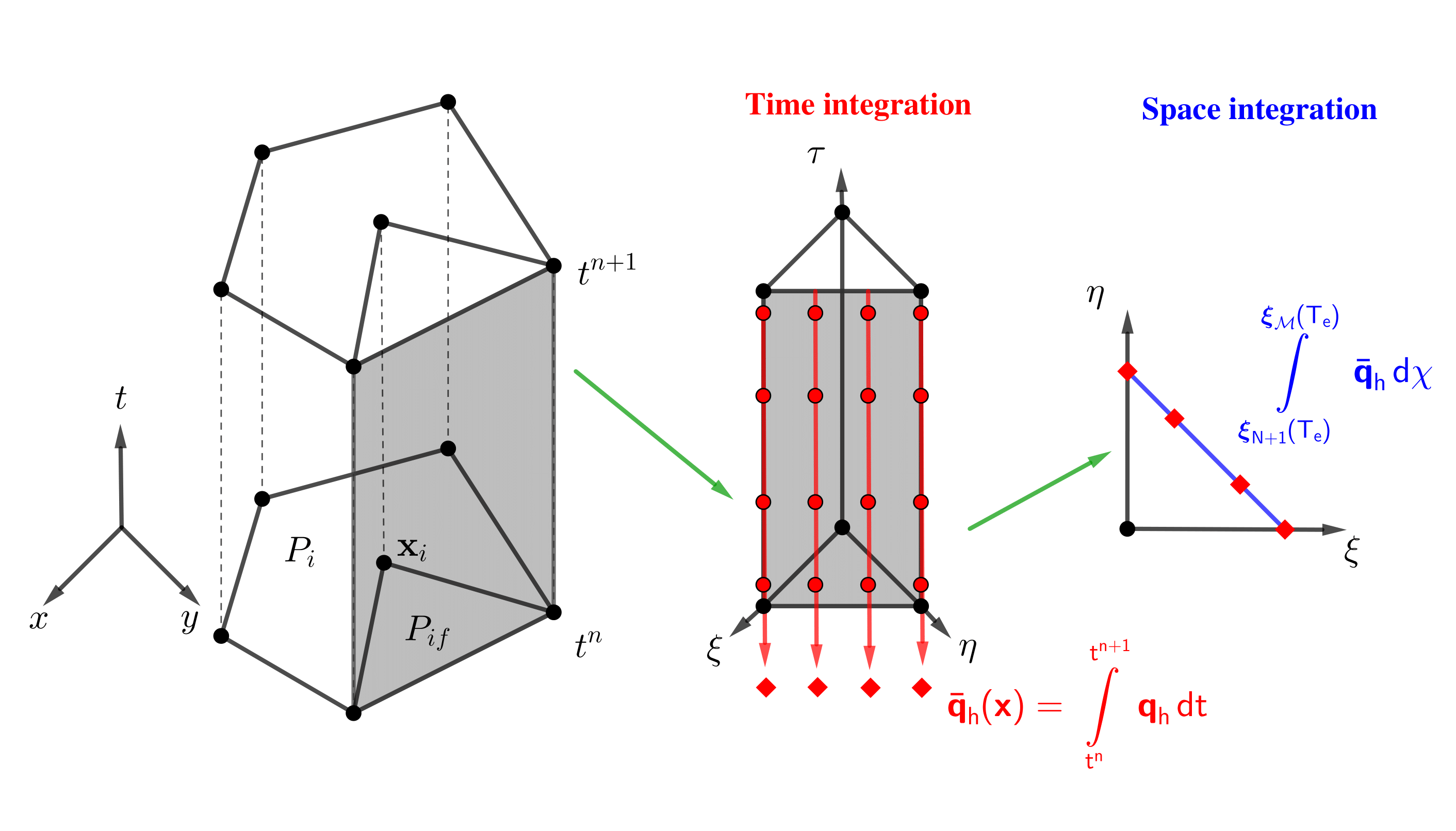} \\ 
		\end{tabular} 
		\caption{The space-time subcell $P_{if}\times [t^n;t^{n+1}]$ (left) is mapped to the reference space-time element $\tilde{T}_e=T_e \times [0;1]$ where the expansion coefficients $\hat{\q}_{\ell,i}$ for $N=3$ are integrated in time, thus yielding $\hat{\bar{\q}}_{\ell,i}$ in the reference element (middle). Next, these time integrated values are used for designing a quadrature-free scheme with flux integration in space of the degrees of freedom $\hat{\bar{\q}}_{\ell,i}$ according to the definition of the spatial integration matrix $\T_{\xxi}$~\eqref{eqn.DGmatrices} (right).}
		\label{fig.space-time_int}
	\end{center}
\end{figure}

Moreover, the definition of the space-time basis functions in terms of a tensor product allows time integration to be performed immediately once system~\eqref{eqn.PDEweak2} is solved and the full space-time predictor solution $\q_h$ is available for each element $P_i$. Indeed, with this tensor product definition, each degree of freedom over the \textit{spatial} reference element $T_e$ has its $N+1$ counterparts in \textit{time}, as depicted in Figure~\ref{fig.space-time_int}. One can therefore compute the time integrated predictor solution $\bar{\q}_h$ \rd{and the time integrated flux $\bar{\F}_h$} as
\begin{equation} 
\begin{aligned} 
& \bar{\q}_h(\x) = \int \limits_{t^n}^{t^{n+1}} \q_h \, dt = \dt \, \sum \limits_{f \in \mathcal{F}_i} \T_{\tau} \, \hat \q_{\ell,if}, \qquad \x \in P_i, \qquad \T_{\tau}:=\int \limits_{0}^{1} \theta_\ell(\tau) \, d\tau, \\
& \rd{ \bar{\F}_h = (\bar{\f}_h, \bar{\g}_h) \ \text{ with } \
\bar{\f}_h(\x) = \int \limits_{t^n}^{t^{n+1}} \f_h^* \, dt = \dt \, \sum \limits_{f \in \mathcal{F}_i} \T_{\tau} \, \hat \f_{\ell,if}^* \ \text{ and } \ 
\bar{\g}_h(\x) = \int \limits_{t^n}^{t^{n+1}} \g_h^* \, dt = \dt \, \sum \limits_{f \in \mathcal{F}_i} \T_{\tau} \, \hat \g_{\ell,if}^*, }
\end{aligned}
\end{equation} 
with $\T_{\tau}$ denoting a universal matrix which is given by the integral of the space-time basis functions over the reference time interval $[0;1]$. Furthermore, notice that the time-integrated predictor solution $\bar{\q}_h(\x)$ does depend only on space and they can thus be expressed using the \textit{same} agglomerated basis functions $\phi_\ell$ adopted for the numerical solution $\u_h$ in~\eqref{eqn.uh}: 
\begin{equation}
\bar{\q}_h(\x) = \sum \limits_{\ell=1}^{\mathcal{N}_i} \phi_\ell(\x) \, \hat{\bar{\q}}_{\ell,i} 
:= \phi_\ell(\x) \, \hat{\bar{\q}}_{\ell,i} , \qquad \x \in \rd{P_i}.
\label{eqn.qhTAv}
\end{equation}

\subsection{Fully discrete quadrature-free one-step ADER DG schemes} \label{ssec.corrector}
The time-integrated predictor solution $\bar{\q}_h$, obtained with the predictor step, does not account for the neighboring flux contributions, so that a \textit{corrector} step must follow, which is based on the design of an explicit quadrature-free one-step DG solver. The governing PDE~\eqref{eqn.PDE} is first multiplied by a test function $\phi_k$ of the same form of the basis function $\phi_\ell$ of~\eqref{eqn.uh}, then it is integrated over the space-time control volume $P_i \times [t^n;t^{n+1}]$, obtaining the following weak problem
\begin{equation}
\int \limits_{t^n}^{t^{n+1}} \int \limits_{P_i} \phi_k \, \left( \frac{\partial \Q}{\partial t} + \nabla \cdot \F(\Q,\nabla \Q) \right) \, d\x \, dt = \mathbf{0}. 
\label{eqn.DG1}
\end{equation}  
As typically done for DG schemes, the flux divergence term is integrated by parts in space, hence obtaining
\begin{equation}
\int \limits_{t^n}^{t^{n+1}} \int \limits_{P_i} \phi_k \, \frac{\partial \Q}{\partial t}  \, d\x \, dt + \int \limits_{t^n}^{t^{n+1}} \int \limits_{\partial P_i} \phi_k \, \F(\Q,\nabla \Q) \cdot \mathbf{n} \, dS \, dt - \int \limits_{t^n}^{t^{n+1}} \int \limits_{P_i} \nabla \phi_k \cdot \F(\Q,\nabla \Q) \, d\x \, dt = \mathbf{0}, 
\label{eqn.DG2}
\end{equation}
with $\partial P_i$ denoting the Voronoi cell boundary and $\mathbf{n}$ being the outward pointing unit normal vector defined on $\partial P_i$. The fluxes in~\eqref{eqn.DG2} are computed using the predictor solution $\bar{\q}_h$, which already accounts for time integration and \textit{not} relying on the purely spatial solution $\u_h$, as done in classical multi-stage DG schemes based on Runge-Kutta time stepping techniques. 

To devise a quadrature-free scheme, the weak form~\eqref{eqn.DG2} must be reformulated over the space-time reference system ~$\tilde{\xxi}$.
To do that, we start by integrating the first term of~\eqref{eqn.DG2} in space and time and introducing the definitions~\eqref{eqn.uh} and~\eqref{eqn.qhTAv}. 
Then, the physical flux at the the cell interface $\partial P_i$ is substituted by the numerical flux function $\mathcal{G}_{\rd{f}}(\bar{\q}_h^+,\bar{\q}_h^-)$, and the remaining terms are rewritten with respect to the space-time reference element $\tilde{T}_e$ for each subcell $P_{if}$. Thus, we get the following expression of the fully discrete quadrature-free one-step DG scheme:
\begin{eqnarray}
\M_i \left( \hat{\mathbf{u}}^{n+1}_{\ell,i} - \hat{\mathbf{u}}^{n}_{\ell,i} \right) + \sum \limits_{f \in \mathcal{F}_i} |\partial P_{if}| \, \T_{\xxi} \, \mathcal{G}_f(\bar{\q}_h^+,\bar{\q}_h^-) = \sum \limits_{f \in \mathcal{F}_i} |J_{if}| \, \left[ \V_{\xi} \f^*(\bar{\q}_h) + \V_{\eta} \g^*(\bar{\q}_h) \right],
\label{eqn.DGQF}
\end{eqnarray}
where the element-wise mass matrix writes
\begin{equation}
\M_i = |P_i| \, \sum \limits_{f\in \mathcal{F}_i} \int \limits_{T_e} \phi_k \phi_\ell \, d\xxi,
\label{eqn.Mi}
\end{equation}
while the remaining universal matrices are given by
\begin{equation}
\T_{\xxi} = \int \limits_{\xxi_{N+1}(T_e)}^{\xxi_{\mathcal{M}}(T_e)} \phi_k \phi_\ell \, d\chi, \qquad \V_{\xi} = \int \limits_{T_e} \frac{\partial \phi_k}{\partial \xi} \phi_\ell \, d\xxi, \qquad \V_{\eta} = \int \limits_{T_e} \frac{\partial \phi_k}{\partial \eta} \phi_\ell \, d\xxi.
\label{eqn.DGmatrices}
\end{equation}
As for the combined temporal stiffness matrix $\K_{1,i}$ in the predictor step, the inverse of the mass matrix $\M_i^{-1}$, which is the \textit{only element-dependent} matrix involved in the corrector step, can be conveniently computed and stored for each element one and for all at the beginning of the algorithm. The flux matrix $\T_{\xxi}$ is \textit{always} defined along the reference edge with the node coordinates $\xxi_{N+1}$ and $\xxi_{\mathcal{M}}$ as extrema (we recall indeed that each subcell $P_{if}$ is always mapped to the reference element $T_e$ with the cell barycenter $\x_i$ as the first local node $\xxi_1$ according to~\eqref{eqn.xietaTransf}, see also Figure~\ref{fig.Te}).

To account for the discontinuities arising in the boundary integral of~\eqref{eqn.DG2}, we rely on the usage of a Riemann solver $\mathcal{G}$ at the element interfaces~\cite{ToroBook}. The right and left data necessary for computing $\mathcal{G}$ are then given by the right and left states $(\bar{\q}_h^+,\bar{\q}_h^-)$ obtained within the predictor step, thus yielding the formal order of accuracy $\mathcal{O}(N+1)$. The left state $\bar{\q}_h^-$ refers to the cell under consideration, i.e. $P_i$, whereas the right state $\bar{\q}_h^+$ is related to the neighbor cell $P_j$ which shares the common face $f$ across which the flux is evaluated. A simple and very robust Rusanov flux~\cite{Rusanov:1961a} is adopted, which is modified to simultaneously include both the convective and the viscous terms~\cite{MunzDiffusionFlux}:
\begin{equation}
\mathcal{G}_f = \frac{1}{2} \left( \rd{\bar{\F}_h^+} + \, \rd{\bar{\F}_h^-} \right) \cdot \n_{if} - \frac{1}{2} \left( |\lambda^{\max}| + 2 \varepsilon |\lambda_v^{\max}| \right) \left( \bar{\q}_h^+ - \bar{\q}_h^- \right), \qquad \varepsilon = \frac{2N+1}{(h_i+h_j)\sqrt{\frac{\pi}{2}}}.
\label{eqn.rusanov}
\end{equation} 
The quantities $|\lambda^{\max}|$ and $|\lambda_v^{\max}|$ are obtained taking the corresponding maximum absolute value of the eigenvalues in~\eqref{eqn.eigenval} between $\bar{\q}_h^+$ and $\bar{\q}_h^-$, then projected in the face normal direction, i.e. along the unit vector $\n_{if}$. As proposed in~\cite{DumbserKaeser07}, these values are then \textit{frozen} for the flux evaluation across the entire face $f$, making the numerical flux~\eqref{eqn.rusanov} a \textit{linear} function of its \textit{four} arguments $(\bar{\q}_h^+,\bar{\q}_h^-,\bar{\F}_h^+,\bar{\F}_h^-)$ and thus allowing to obtain a \textit{quadrature-free} formulation on general polygonal grids.

\subsection{DG limiter with artificial viscosity} \label{ssec.limiter}
In the case of shock waves and other discontinuities that might occur even when starting the simulation with smooth initial data, because of the non-linearity of the governing PDE system, a limiter technique must be introduced in the DG scheme~\eqref{eqn.DGQF} in order to avoid spurious oscillations and dangerous instabilities.	
Among the variety of techniques presented in literature \cite{Kuzmin_Limiter2010,SonntagDG,DGLimiter1,DGLimiter2,ALEDG,DGCWENO,DeLaRosaMunzDGMHD,Kuzmin_Limiter2020,Shu_Limiter2020,Gaburro_LimiterPNPM,Gassner_Limiter2021}	
that range from classical \textit{a priori} limiters to novel \textit{a posteriori} limiting techniques, and that exploit the combination of different schemes, bounds preserving approaches and/or different level of refinements,  
here we have chosen to introduce a simple artificial viscosity method, inspired from \cite{PerssonAV,GassnerMunzAV,TavelliCNS,Hesthaven_LimiterAV2011,Bassi_LimiterAV2018,Bassi_LimiterAV2020}, 
that we apply only to those cells that need limitation, i.e. those which are detected as \quotew{troubled}. 
Therefore, our limiter procedure is made of two steps: i) detection and ii) limiting.

To detect these troubled elements, the flattener variable $\beta$ proposed in~\cite{BalsaraFlattener} is used as indicator at the beginning of each time step. This detector works by comparing the divergence of the velocity field $\nabla \cdot \v$ in the element $P_i$ with the minimum of the sound speed $c=\sqrt{\gamma R T}$ among the element $P_i$ itself and its neighborhood. Let us introduce an estimation of the divergence of the velocity field given by
\begin{equation}
\nabla \cdot \v|_i = \frac{1}{|P_i|}\sum \limits_{f \in \mathcal{F}_i}{ |\partial P_{if}|\left(\v^+ - \v^- \right) \cdot \n_{if} }, \qquad c_{s,\min} = \min \limits_{f \in \mathcal{F}_i}{\left(c_{s}^-,c_{s}^+\right)},
\label{eqn.divV}
\end{equation}
with the superscripts $(+,-)$ denoting the right and left quantities across the face $f$. The flattener variable $\beta_i$, $i=1,\ldots,N_E$ is now computed according to~\cite{BalsaraFlattener} as
\begin{equation}
\beta_i = \min {\left[ 1, \max {\left(0, -\frac{\nabla \cdot \v_i + m_1 c_{s,\min}}{m_1 c_{s,\min}}\right)}\right]},
\label{eqn.flattener}
\end{equation}
with the coefficient $m_1$ that is set to the value $m_1=0.1$ as done in~\cite{BalsaraFlattener}.

If the detector is activated on a cell $P_i$ then this cell is said to be \quotew{troubled} and some artificial viscosity $\mu_{add,i}$ is added to the physical viscosity $\mu$, hence resulting in the effective viscosity $\mu_i=\mu_{add,i}+\mu$ used in the Navier-Stokes equations~\eqref{eqn.NSEterms}. The additional viscosity
$\mu_{add,i}$ is determined so that a resulting unity mesh Reynolds number $Re_i$ is assigned to the troubled cells, that is
\begin{equation}
Re_i = \frac{\rho |\lambda^{\max}_i| h_s}{\mu_i},
\end{equation}
with $|\lambda^{\max}_i|$ given by the maximum absolute value of the convective eigenvalues in~\eqref{eqn.eigenval} and $h_s=h_i/(2N+1)$ representing a rescaled mesh spacing according to~\cite{TavelliCNS}. Consistently, the corresponding artificial heat conduction coefficient $\kappa_i$ in the heat flux is obtained by setting the Prandtl number to $Pr_i = \mu_i \gamma c_v/\kappa_i = 1$.

\section{Numerical results} \label{sec.test}
The aim of this section is to describe and show the numerical results for a series of test cases solved using our novel DG schemes based on agglomerated subgrid continuous finite element basis functions. The abbreviation ADER DG-AFE (Agglomerated Finite Element) is used for referring to the method presented in this work, while the label ADER DG-M (Modal) refers to the numerical method forwarded in~\cite{ArepoTN} and based on standard modal Taylor-type basis functions written in terms of the physical coordinates. The CFL number in~\eqref{eqn.timestep} is set to $\CFL=0.5$ for all simulations and the obtained numerical
solutions are compared against exact or numerical reference solutions available in the literature. If not otherwise specified, the ratio of specific heats is $\gamma=1.4$ and the gas constant is $R=1$, thus the specific heat capacity at constant volume is set to $c_v=2.5$. Furthermore, the limiting strategy presented in Section~\ref{ssec.limiter} is not activated by default. The initial condition of the flow field is typically given in terms of the vector of primitive variables $\mathbf{P}(\x,t)=(\rho,u,v,p)$ and all simulations are run on 64 CPUs with MPI parallelization.

Table~\ref{tab.efficiency} reports a comparison in terms of computational efforts between the AFE and the M version of the ADER DG method on Voronoi tessellations for all the test problems proposed in this section. For the ADER DG-AFE scheme the time required for an element update within one time step can be up to $\approx 2.8$ times faster than the corresponding ADER DG-M algorithm, thus a remarkable improvement in the overall performance of the novel schemes is achieved, while maintaining very accurate resolution properties as demonstrated in the following test suite.

\begin{table}[!t]  
	\caption{Profiling of both ADER DG-AFE and ADER DG-M schemes solving all test cases proposed in this work (EP=Explosion Problem, SP=Stokes Problem, VS=Viscous Shock, TG=Taylor-Green, ML=Mixing Layer). The order of accuracy of the scheme is $\mathcal{O}(N)$, the number of elements involved in the spatial discretization is denoted with $N_E$, the number of time steps needed to reach the final time is given by $N_{\Delta t}$, while the absolute time of each simulation measured in seconds $[s]$ is provided for both schemes in the columns $T_{AFE}$ and $T_{M}$. The time needed for the update of one element within one time step is indicated with $\tau_{AFE}$ and $\tau_{M}$. Finally, the efficiency ratio is evaluated as $\mathcal{R}_{E}=\tau_{M}/\tau_{AFE}$.}  
	\begin{center} 
		\begin{small}
			\renewcommand{\arraystretch}{1.1}
			\begin{tabular}{lcr|rcc|rcc|c}
				\multirow{2}{*}{Test} & \multirow{2}{*}{$\mathcal{O}(N+1)$} & \multirow{2}{*}{$N_E$} & \multicolumn{3}{c|}{ADER DG-AFE} & \multicolumn{3}{c|}{ADER DG-M} & \\
				& & & $N_{\Delta t}$ & $T_{AFE}$ & $\tau_{AFE}$ & $N_{\Delta t}$ & $T_{M}$ & $\tau_{M}$ & $\mathcal{R}_{E}$\\ 
				\hline
				EP  & $\mathcal{O}(3)$ & 25648 & 3856 & 2.237E+05 & 2.262E-03 & 3860 & 3.088E+05 & 3.119E-03 & 1.38 \\
				SP ($\mu=10^{-2}$)  & $\mathcal{O}(4)$ & 358 & 90093 & 4.388E+05 & 1.361E-02 & 90098 & 1.219E+06 & 3.780E-02 & 2.78 \\
				SP ($\mu=10^{-3}$)  & $\mathcal{O}(4)$ & 358 & 11546 & 6.762E+04 & 1.636E-02 & 11546 & 1.792E+05 & 4.366E-02 & 2.65 \\ 
				SP ($\mu=10^{-4}$)  & $\mathcal{O}(4)$ & 358 & 3689 & 2.694E+04 & 2.040E-02 & 3689 & 6.172E+04 & 4.674E-02 & 2.29 \\
				VS  & $\mathcal{O}(4)$ & 1120 & 73372 & 7.129E+05 & 8.675E-03 & 73372 & 1.746E+06 & 2.125E-02 & 2.45 \\
				TG  & $\mathcal{O}(4)$ & 2916 & 8272 & 2.501E+05 & 1.037E-02 & 8272 & 6.486E+05 & 2.689E-02 & 2.59 \\
				ML  & $\mathcal{O}(3)$ & 15723 & 65786 & 1.759E+06 & 1.700E-03 & 65787 & 3.327E+06 & 3.216E-03 & 1.89 \\
			\end{tabular}
		\end{small}
	\end{center}
	\label{tab.efficiency}
\end{table}

\subsection{Numerical convergence studies} \label{ssec.ShuVortex}
\begin{figure}[!htbp]
	\begin{center}
		\begin{tabular}{ccc}  
			\includegraphics[width=0.33\textwidth]{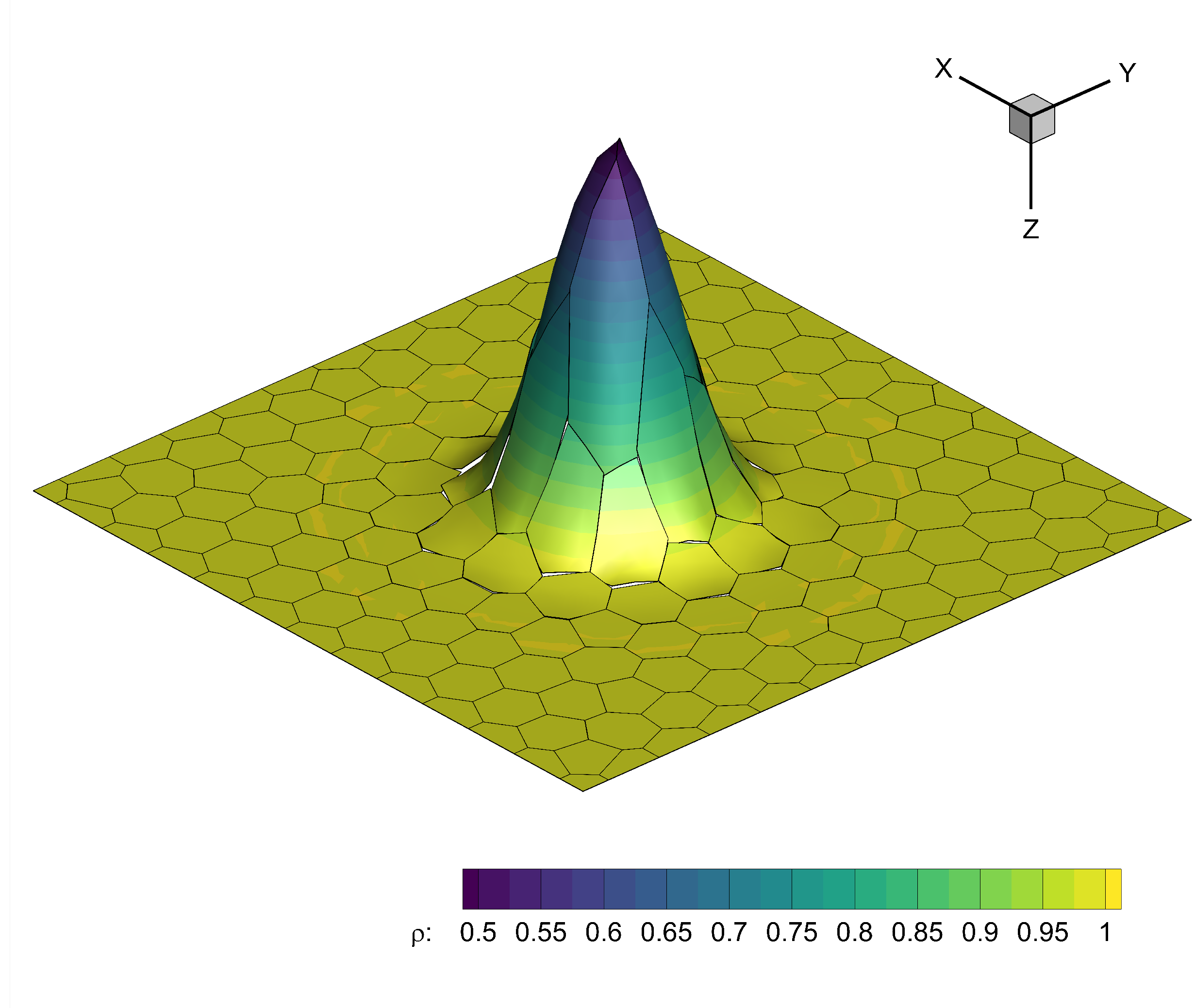} &
			\includegraphics[width=0.33\textwidth]{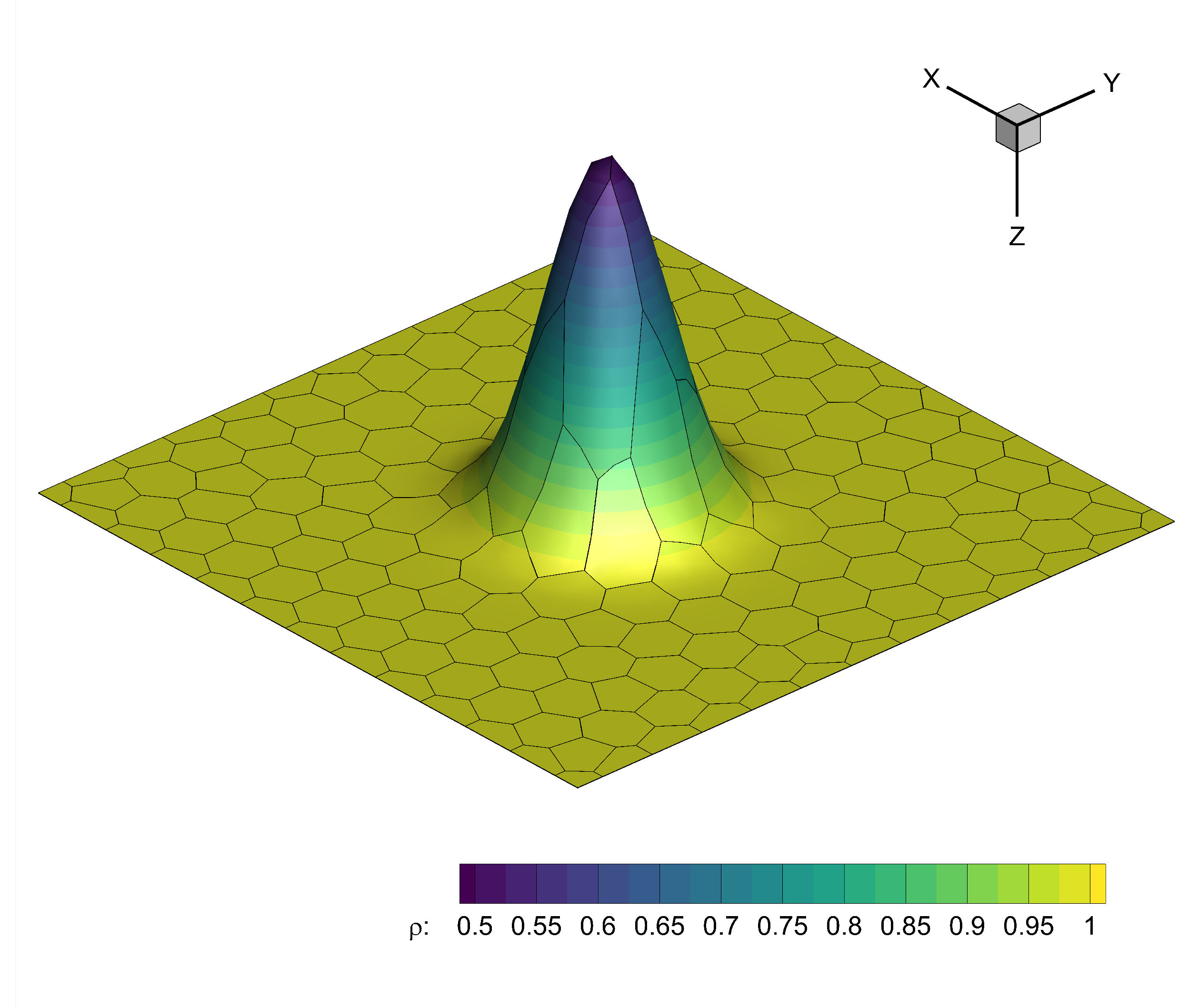} &
			\includegraphics[width=0.33\textwidth]{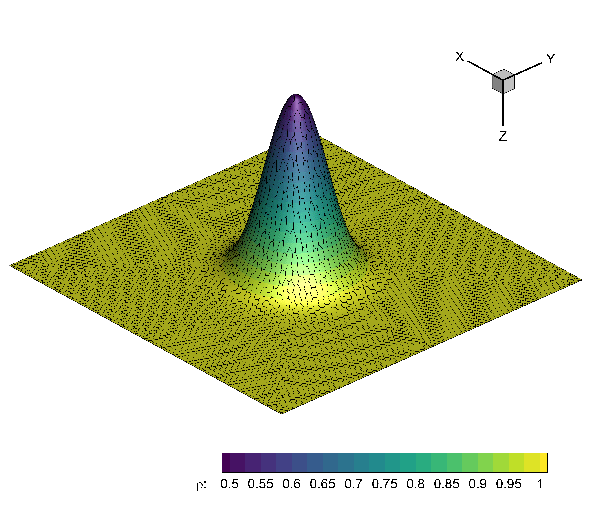} \\ 
		\end{tabular} 
		\caption{$L_2$ projection of the density distribution of the initial condition for the isentropic vortex test case with $N=2$ on a coarse mesh with characteristic size $h=10/12$. Left: agglomerated finite element basis functions. Middle: modal basis functions. Right: reference solution obtained on a computational grid with mesh spacing $h=10/64$.}
		\label{fig.L2IC}
	\end{center}
\end{figure}

To study the numerical convergence of the novel ADER DG-AFE schemes, an isotropic vortex problem is considered according to the setup initially proposed in~\cite{HuShuTri}. This problem
has an exact and smooth solution for the compressible Euler equations, thus viscosity and heat conduction coefficients are neglected, i.e. $\mu=\kappa=0$. The computational domain is the square $\Omega=[0;10]^2$ with periodic boundaries and the general radial coordinate is given by $r=\sqrt{\x^2}$. The vector of primitive variables at the initial time is defined as
\begin{equation}
\label{ShuVortIC}
\mathbf{P}(\x,0) = (1+\delta \rho, 1+\delta u, 1+\delta v, 1+\delta p),
\end{equation}  
where the perturbations for pressure $\delta p$, density $\delta \rho$ and velocity $(\delta u, \delta u)$ are
\begin{equation}
\label{ShuVortDelta}
\delta p = (1+\delta T)^{\frac{\gamma}{\gamma-1}}-1,\quad \delta \rho = (1+\delta T)^{\frac{1}{\gamma-1}}-1,\quad 
\left(\begin{array}{c} \delta u \\ \delta v \end{array}\right) = \frac{\epsilon}{2\pi}e^{\frac{1-r^2}{2}} \left(\begin{array}{c} -(y-5) \\ \phantom{-}(x-5) \end{array}\right), 
\end{equation}
with $\delta T$ being the temperature perturbation that reads
\begin{equation}
\delta T    = -\frac{(\gamma-1)\epsilon^2}{8\gamma\pi^2}e^{1-r^2}.
\end{equation}
The vortex strength is set to $\epsilon=5$, while the perturbation of entropy $S=p/\rho^{\gamma}$ is assumed to be zero. The initial density distribution is shown in Figure~\ref{fig.L2IC} on a very coarse mesh for ADER DG-AFE and ADER DG-M schemes, and a comparison against the reference solution on a very fine grid is proposed, highlighting the benefits in terms of resolution that are achieved by the novel AFE approach. Moreover, the errors related to the $L_2$ projection of the initial condition~\eqref{ShuVortDelta} are reported in Table~\ref{tab.Linf_IC}, which demonstrates that, in $L_{\infty}$ norm, the DG-AFE schemes are up to one order of magnitude more accurate than the corresponding DG-M version. 

\begin{table}[!t]  
	\caption{Errors related to the $L_2$ projection of the initial condition for the isentropic vortex test case measured in $L_{\infty}$ norm for $\rho$ (density), $u$ (horizontal velocity) and $p$ (pressure). A coarse mesh of size $h=10/12$ is used and the errors are reported for third and fourth order DG schemes with AFE and Modal basis functions.}  
	\begin{center} 
		\begin{small}
			\renewcommand{\arraystretch}{1.0}
			\begin{tabular}{c|cc|cc}
				& \multicolumn{2}{c|}{$\mathcal{O}(3)$} & \multicolumn{2}{c}{$\mathcal{O}(4)$} \\
				\hline
				\hline
				Variable & ADER DG-AFE & ADER DG-M & ADER DG-AFE & ADER DG-M \\
				\hline
				$\rho$ & 2.818E-03 & 1.058E-02 & 2.836E-04 & 2.053E-03 \\
				$u$    & 6.511E-03 & 2.817E-02 & 8.510E-04 & 7.678E-03 \\
				$p$    & 2.332E-02 & 8.514E-02 & 2.151E-03 & 2.075E-02 \\
			\end{tabular}
		\end{small}
	\end{center}
	\label{tab.Linf_IC}
\end{table}

\begin{figure}[!htbp]
	\begin{center}
		\begin{tabular}{cc}  
			\includegraphics[width=0.47\textwidth]{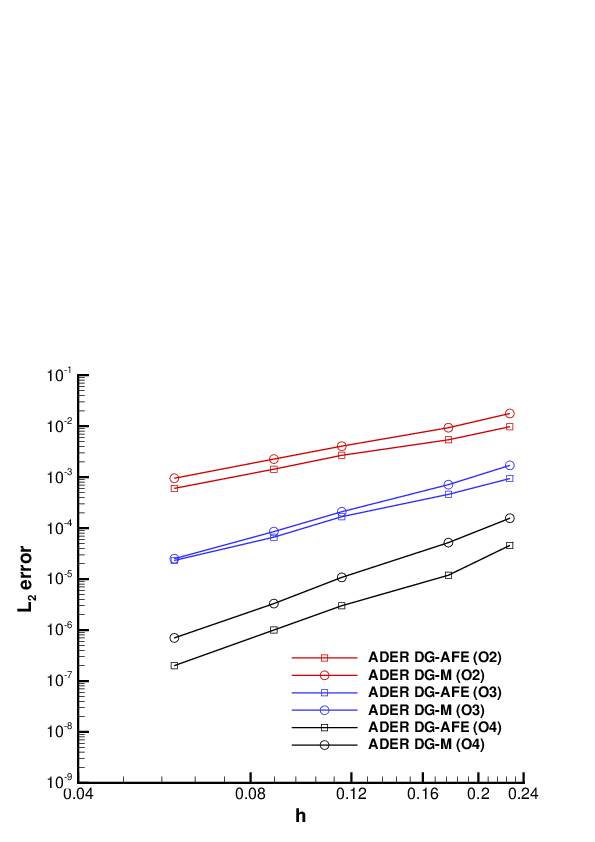} &
			\includegraphics[width=0.47\textwidth]{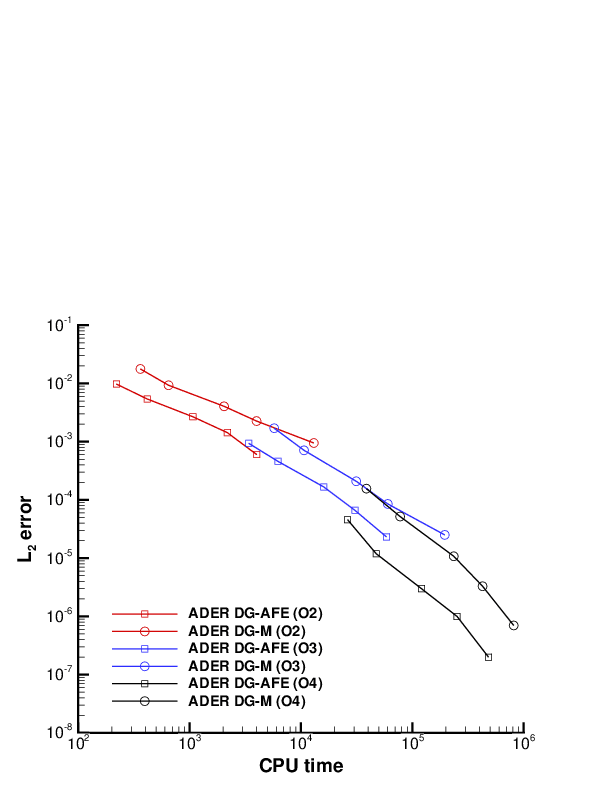} \\ 
		\end{tabular} 
		\caption{\rd{Comparison between ADER DG-AFE and ADER DG-M schemes from second up to fourth order of accuracy. Left: dependency of the error norm on the mesh size. Right: dependency of the error norm on the CPU time.}}
		\label{fig.L2-h-time}
	\end{center}
\end{figure}

The exact solution $\Q_e(\x,t_f)$ at the final time of the simulation $t_f=1$ is simply given by the time-shifted initial condition, thus $\Q_e(\x,t_f) = \Q(\x-\v_c \, t_f,0)$, with the convective velocity of the vortex $\v_c=(u_c,v_c)=(1,1)$ according to the initial condition~\eqref{ShuVortIC}. The computational domain is discretized with a sequence of refined unstructured Voronoi meshes of characteristic size $h(\Omega)=\max_i h_i$ and the corresponding error is measured in $L_2$  norm as 
\begin{equation}
\varepsilon_{L_2} = \sqrt{ \int \limits_{\Omega} \left( \Q_e(\x,t_f) - \u_h(\x,t_f) \right)^2 d\x }, 
\label{eqnL2error}
\end{equation} 
with $\u_h(\x,t_f)$ representing the numerical solution of the DG scheme at the final time. Table~\ref{tab.convRates} shows the convergence results up to fourth order of accuracy in space and time for both ADER DG-AFE and ADER DG-M schemes, with the associated computational times. The novel algorithm systematically provides lower errors with more computational efficiency, which is also highlighted by the error versus CPU time plots depicted in Figure~\ref{fig.L2-h-time}. 

\begin{table}[!t]  
	\caption{Numerical convergence results for the compressible Euler equations using both ADER DG-AFE and ADER DG-M schemes from second up to fourth order of accuracy in space and time. The errors are measured in the $L_2$ norm and refer to the variable $\rho$ (density) at time $t_{f}=1$. The absolute CPU time of each simulation is also reported in seconds $[s]$.}  
	\begin{center} 
		\begin{small}
			\renewcommand{\arraystretch}{1.2}
			\begin{tabular}{c|ccc|ccc}
				\multicolumn{1}{c}{} & \multicolumn{3}{|c}{ADER DG-AFE} & \multicolumn{3}{|c}{ADER DG-M} \\
				\hline
				$h(\Omega)$ & $\rho_{L_2}$ & $\mathcal{O}(\rho_{L_2})$ & CPU time & $\rho_{L_2}$ & $\mathcal{O}(\rho_{L_2})$ & CPU time \\ 
				\hline
				\hline
				& \multicolumn{6}{|c}{Order of accuracy: $\mathcal{O}(2)$} \\
				2.270E-01 & 9.758E-03 & -   & 2.203E+02 & 1.775E-02 & -   & 3.616E+02 \\
				1.773E-01 & 5.406E-03 & 2.4 & 4.174E+02 & 9.322E-03 & 2.6 & 6.472E+02 \\
				1.155E-01 & 2.680E-03 & 1.6 & 1.070E+03 & 4.055E-03 & 1.9 & 2.039E+03 \\
				8.786E-02 & 1.431E-03 & 2.3 & 2.182E+03 & 2.262E-03 & 2.1 & 3.989E+03 \\
				\rd{5.888E-02} & \rd{6.014E-04} & \rd{2.2} & \rd{4.008E+03} &
				\rd{9.504E-04} & \rd{2.2} & \rd{1.304E+04} \\
				& \multicolumn{6}{c}{Order of accuracy: $\mathcal{O}(3)$} \\
				2.270E-01 & 9.369E-04 & -   & 3.395E+03 & 1.704E-03 & -   & 5.750E+03 \\
				1.773E-01 & 4.598E-04 & 2.9 & 6.214E+03 & 7.121E-04 & 3.5 & 1.065E+04 \\
				1.155E-01 & 1.674E-04 & 2.4 & 1.602E+04 & 2.095E-04 & 2.9 & 3.133E+04 \\
				8.786E-02 & 6.637E-05 & 3.4 & 3.049E+04 & 8.542E-05 & 3.3 & 6.020E+04 \\
				\rd{5.888E-02} & \rd{2.317E-05} & \rd{2.6} & \rd{5.853E+04} &
				\rd{2.509E-05} & \rd{3.1} & \rd{1.954E+05} \\
				& \multicolumn{6}{c}{Order of accuracy: $\mathcal{O}(4)$} \\
				2.270E-01 & 4.578E-05 & -   & 2.621E+04 & 1.563E-04 & -   & 3.868E+04 \\
				1.773E-01 & 1.194E-05 & 5.4 & 4.745E+04 & 5.195E-05 & 4.5 & 7.766E+04 \\
				1.155E-01 & 2.963E-06 & 3.3 & 1.207E+05 & 1.085E-05 & 3.7 & 2.354E+05 \\
				8.786E-02 & 1.039E-06 & 3.8 & 2.522E+05 & 3.332E-06 & 4.3 & 4.270E+05 \\
				\rd{5.888E-02} & \rd{2.303E-07} & \rd{3.8} & \rd{4.828E+05} &
				\rd{6.609E-07} & \rd{4.0} & \rd{8.128E+05} \\
			\end{tabular}
		\end{small}
	\end{center}
	\label{tab.convRates}
\end{table}

\gf{The usage of polygonal control volumes allows for larger time steps because the characteristic mesh size is typically larger than the one corresponding to the associated primary triangular mesh. Indeed, Table \ref{tab.hPoly_VS_hTri} reports the mesh spacing for some of the computational grids used to study the numerical convergence of the novel schemes, confirming that the Voronoi polygons exhibit a larger incircle diameter, according to the definition \eqref{eqn.hi}, compared to the triangles of the underlying Delaunay triangulation. Furthermore, the total number of Voronoi polygons is about two times smaller than the number of primal Delaunay triangles. As such, the use of unstructured polygonal meshes is computationally more advantageous compared to classical triangular simplex meshes.} \Rall{Furthermore, Voronoi tessellations demonstrate to be more suitable and robust to be adopted in the context of Arbitrary Lagrangian-Eulerian schemes as shown in \cite{Springel} and \cite{ArepoTN}, where applications to very long and complex astrophysical simulations are proposed, which are characterized by strong differential rotations phenomena.}

\begin{table}[!t]  
	\caption{\gf{Comparison of the average incircle diameter given by \eqref{eqn.hi} between the primary triangular grid ($h^T$) and the corresponding dual Voronoi tessellation ($h^P$) for a sequence of meshes used in the numerical convergence studies. The characteristic mesh size is measured in terms of the number $N_I$ of cells which split each edge of the computational domain $\Omega=[0;10]^2$.}}  
	\begin{center} 
		\begin{small}
			\renewcommand{\arraystretch}{1.0}
			\begin{tabular}{c|cc}
				$N_I$ & $h^P$ & $h^T$  \\
				\hline
				12 & 4.521E-01 & 2.535E-01 \\
				24 & 2.014E-01 & 1.037E-01 \\
				32 & 1.589E-01 & 9.175E-02 \\
			\end{tabular}
		\end{small}
	\end{center}
	\label{tab.hPoly_VS_hTri}
\end{table}

\subsection{Circular explosion problem} \label{ssec.EP2D}
\begin{figure}[!htbp]
	\begin{center}
		\begin{tabular}{cc} 
			\includegraphics[trim= 1 1 1 1, clip, width=0.47\textwidth]{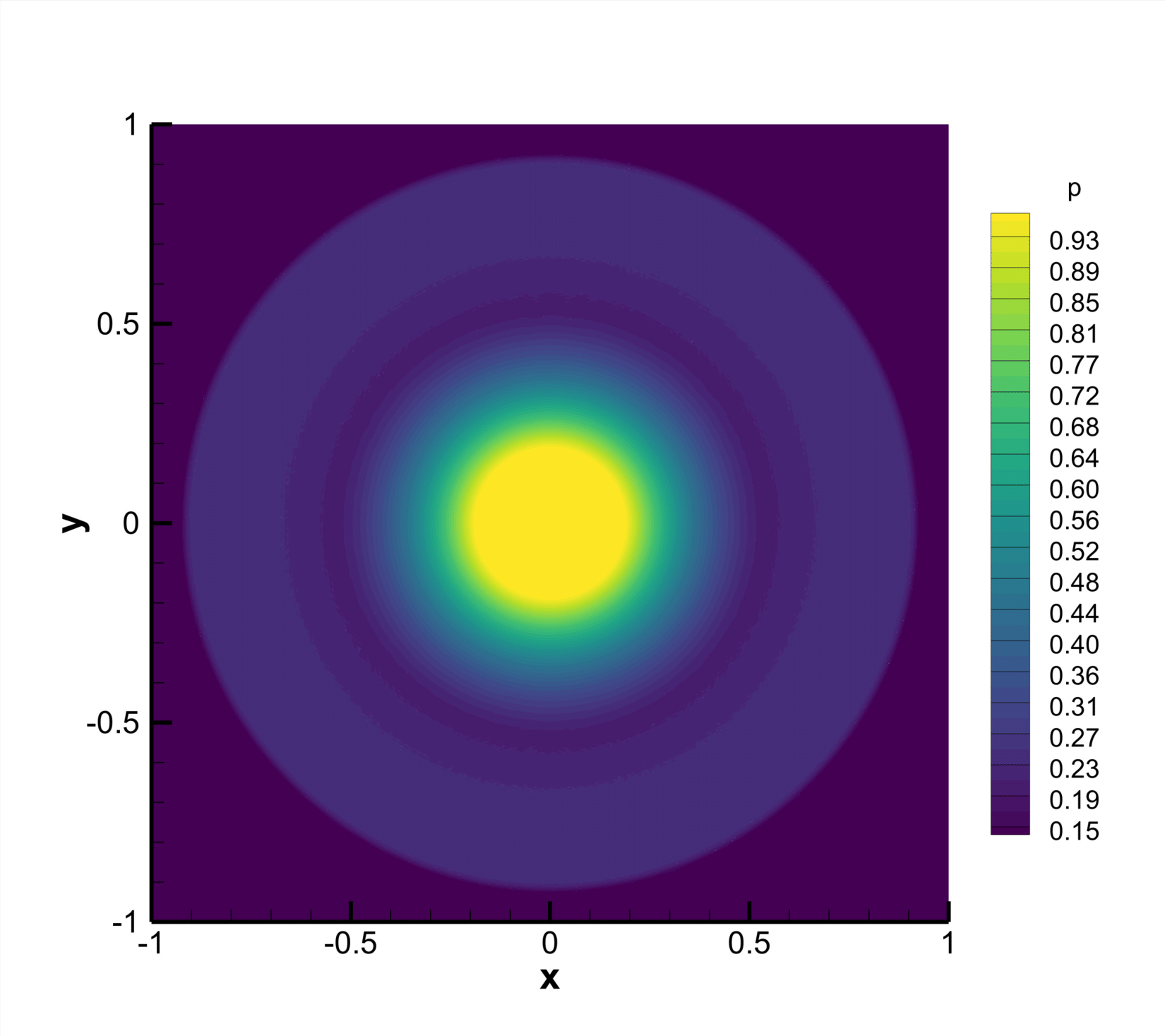}&
			\includegraphics[width=0.47\textwidth]{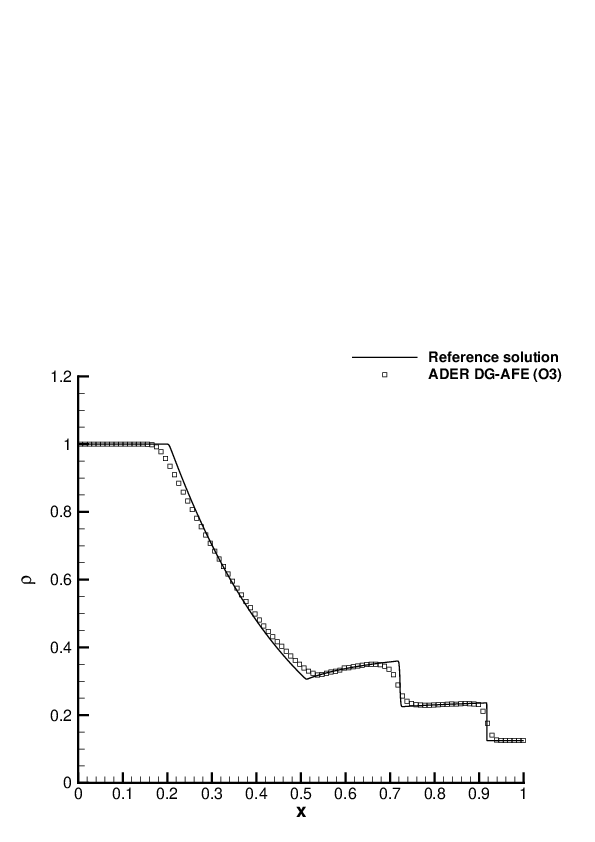} \\				
			\includegraphics[width=0.47\textwidth]{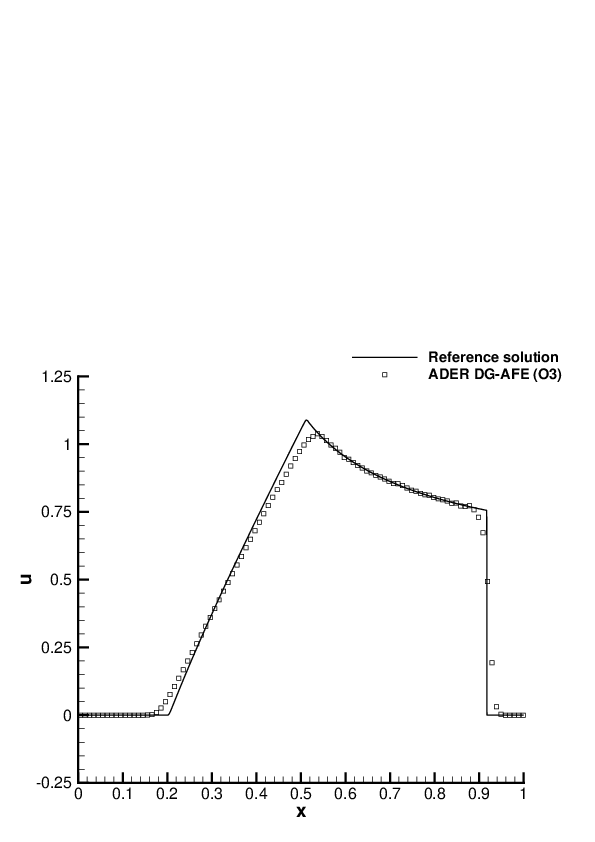} & 
			\includegraphics[width=0.47\textwidth]{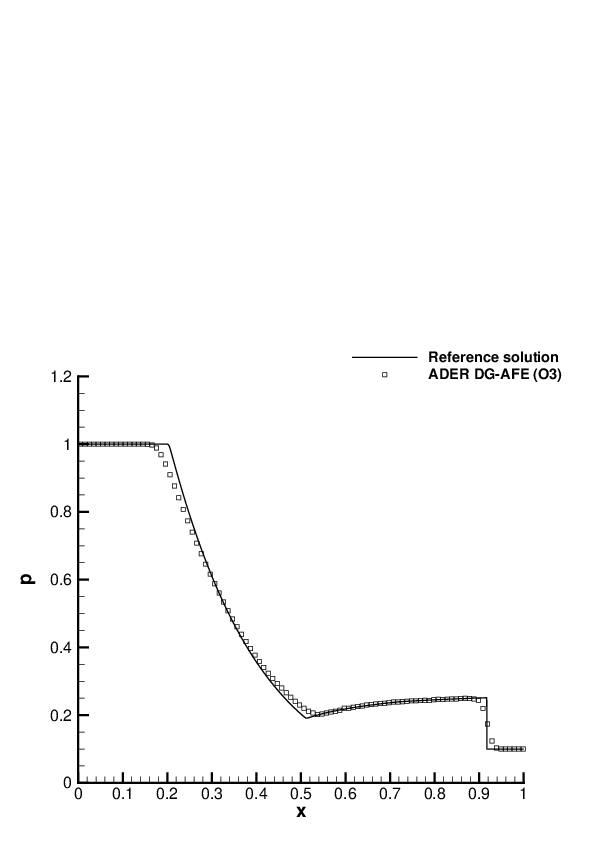} \\
		\end{tabular} 
		\caption{Circular explosion problem at time $t_f=0.25$. Third order numerical results with ADER DG-AFE scheme for density, horizontal velocity and pressure compared against the reference solution
			extracted with a one-dimensional cut of 200 equidistant points along the $x-$direction at $y=0$. Pressure distribution with 40 contour levels in the interval $[0.1;1]$ is shown in the top left panel.}
		\label{fig.EP2D}
	\end{center}
\end{figure}

\begin{figure}[!htbp]
	\begin{center}
		\begin{tabular}{cc}  
			\includegraphics[width=0.47\textwidth]{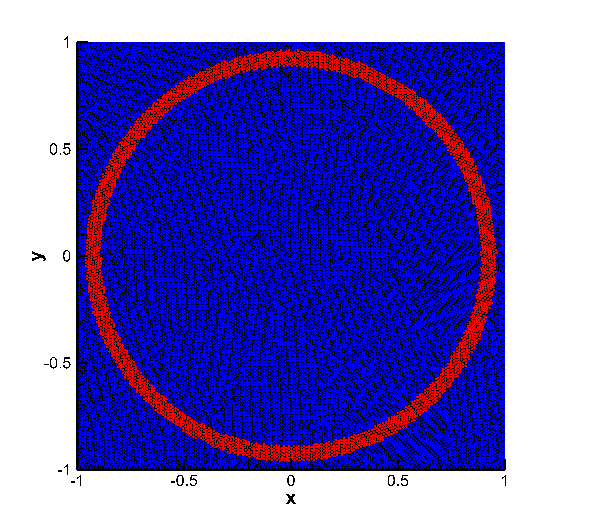} & 
			\includegraphics[width=0.47\textwidth]{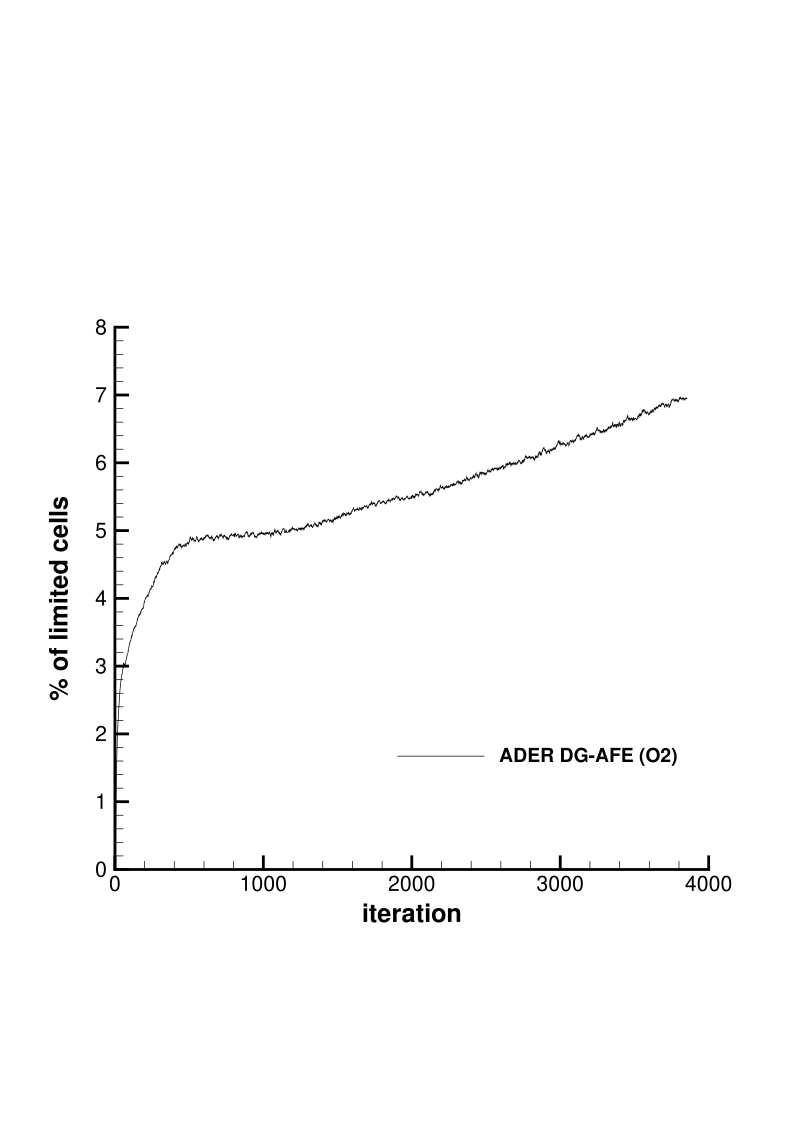} \\ 
		\end{tabular} 
		\caption{Circular explosion problem with third order ADER DG-AFE scheme. Troubled cell map (left) with the limited cells highlighted in red and the unlimited cells colored in blue; percentage of limited cells at each time step (right).}
		\label{fig.EP2Dlim}
	\end{center}
\end{figure}

The second test problem concerns again an inviscid fluid without heat conduction and it consists of a genuinely multidimensional problem which is characterized by a radially symmetric solution. This test case is known as explosion problem and it deals with a cylindrical Riemann problem involving three different types of waves, namely an outward traveling shock front, an inward moving rarefaction fan and a contact wave in between. Due to the presence of a shock wave in the inviscid case, the limiting strategy is activated for this test case. The computational domain is the square $\Omega=[-1;1]^2$ that is assigned transmissive boundary conditions and is meshed with a total number of $N_E=25648$ Voronoi elements. The initial condition is  given in terms of two different states, separated by the circle of radius $R=0.5$: 
\begin{equation}
\mathbf{P}(\x,0) = \left\{ \begin{array}{clcc} \mathbf{P}_i = & (1.0, 0.0, 0.0, 1.0)    & \textnormal{ if } & r \leq R, \\ 
\mathbf{P}_o = & (0.125, 0.0, 0.0, 0.1)  & \textnormal{ if } & r > R,        
\end{array}  \right. 
\end{equation}
where $\mathbf{P}_i$ and $\mathbf{P}_o$ represent the inner and the outer state, respectively, whereas $r=\sqrt{\x^2}$ is the general radial position. In order to avoid nonphysical oscillations at the initial time, where the discontinuous initial data are projected onto the piecewise polynomial approximation space via classical $L_2$ projection, the initial condition is slightly smoothed according to~\cite{TavelliCNS} as follows:
\begin{equation}
\mathbf{P}(\x,0) = \frac{1}{2} \left(\mathbf{P}_o+\mathbf{P}_i\right) + \frac{1}{2} \left(\mathbf{P}_o-\mathbf{P}_i\right) \textnormal{erf} \left( \frac{r-R}{\alpha_0} \right), \qquad \alpha_0=2 \cdot 10^{-2}.
\end{equation}
Due to the cylindrical symmetry of the problem the solution can be compared with an equivalent one-dimensional problem in radial direction $r$ with a geometric source term, see~\cite{ToroBook}. The reference solution is computed at the final time $t_f=0.25$ by solving the compressible Euler equations using a second order MUSCL scheme with the Rusanov flux on a very fine one-dimensional mesh of 15000 points in the radial interval $r\in[0;1]$. Figure~\ref{fig.EP2D} shows the comparison between the reference solution and the numerical solution obtained with the third order version of the ADER DG-AFE scheme as well as a two-dimensional view of the pressure distribution which exhibits excellent symmetry preservation despite the unstructured nature of the computational mesh. A very good agreement can be noticed, and the limiter only acts on those cells which are located at the shock front, that are less than $7\%$ of the total number of elements, as depicted in Figure~\ref{fig.EP2Dlim}.

\subsection{First problem of Stokes} \label{ssec.Stokes}
\begin{figure}[!htbp]
	\begin{center}
		\begin{tabular}{ccc} 
			\includegraphics[width=0.33\textwidth]{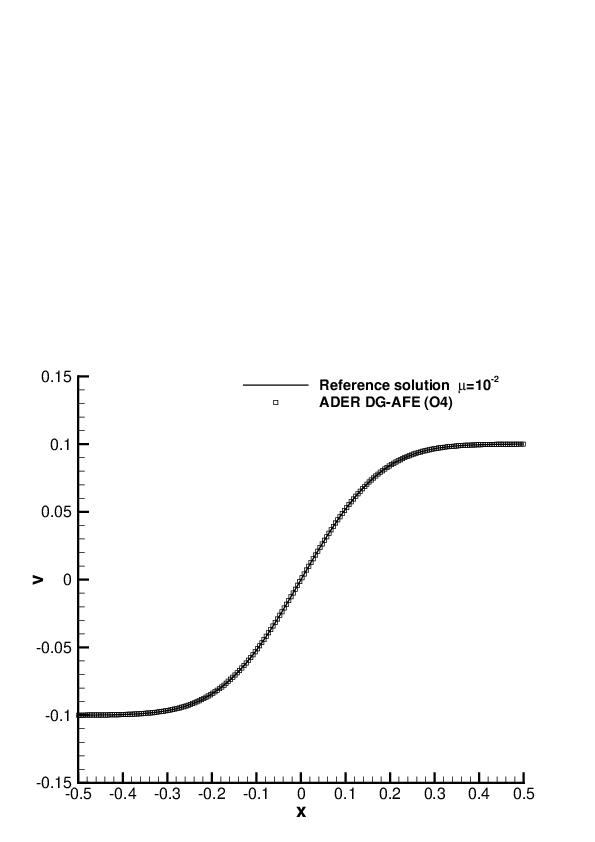} & 
			\includegraphics[width=0.33\textwidth]{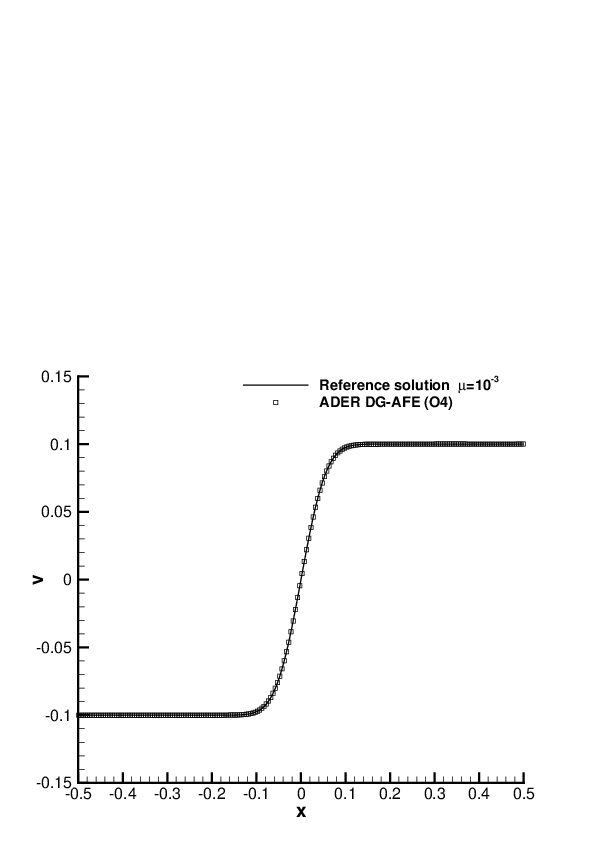} & 
			\includegraphics[width=0.33\textwidth]{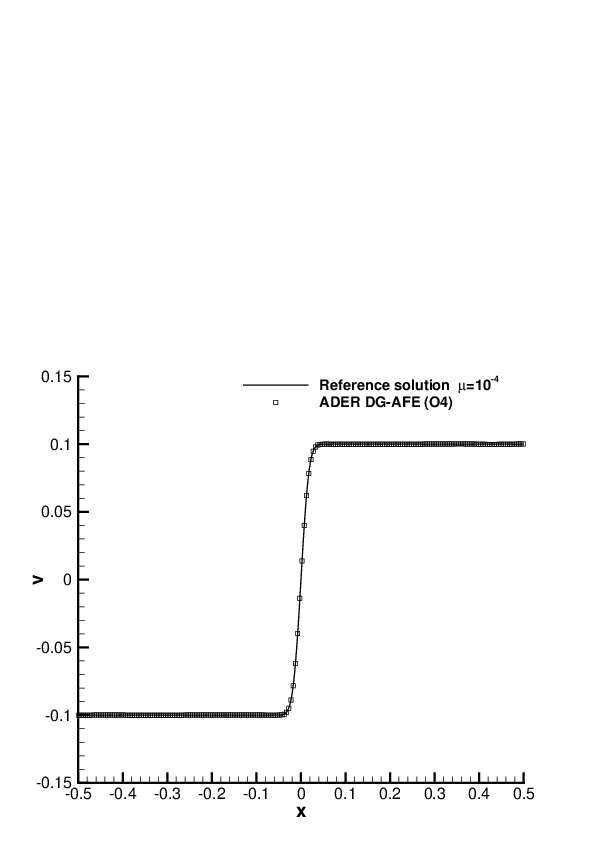} \\
		\end{tabular} 
		\caption{First problem of Stokes at time $t_f=1$. Fourth order numerical results for the vertical component of the velocity obtained with our ADER DG-AFE scheme and compared against the reference solution by extracting a one-dimensional cut of 200 equidistant points along the $x-$direction at $y=0$. Viscosity $\mu=10^{-2}$ (left), $\mu=10^{-3}$ (middle) and $\mu=10^{-4}$ (right).}
		\label{fig.Stokes}
	\end{center}
\end{figure}

We now consider a simple test problem dominated by viscosity effects, namely the first problem of Stokes~\cite{Schlichting}, for which an exact analytical solution of the unsteady Navier-Stokes equations is available. This test describes the time-evolution of an infinite incompressible shear layer, thus the simulation is run at a low Mach number of $M = 0.1$ to obtain an almost incompressible behavior. The computational domain is the channel $\Omega=[-0.5,0.5]\times[-0.05,0.05]$ and is discretized with very few control volumes, namely $N_E=358$. Periodic boundary conditions are prescribed in $y-$direction, while in $x-$direction we enforce the initial condition given by 
\begin{equation} 
\rho = 1, \quad u=0, \quad v=\left\{ \begin{array}{rl}
v_0 & x \leq 0 \\ -v_0 & x > 0 
\end{array} \right. , \quad p=\frac{1}{\gamma}, \qquad v_0 = 0.1.
\end{equation}
Heat conduction is neglected ($\kappa=0$) and the simulation stops at time $t_f=1$. The exact solution of the incompressible Navier-Stokes equations for the velocity component $v$ can be computed as
\begin{equation}
\label{eqn.StokesExact}
v(x,t) = v_0 \, \textnormal{erf} \left(\frac{1}{2} \frac{x}{\sqrt{\mu \, t}}\right).
\end{equation}
This test case is run using the fourth order ADER DG-AFE schemes for three different values of viscosity, namely $\mu=10^{-2}$, $\mu=10^{-3}$ and $\mu=10^{-4}$. The comparison between the reference solution~\eqref{eqn.StokesExact} and the numerical results is presented in Figure~\ref{fig.Stokes}, where one can appreciate an excellent matching between the two for all the considered viscosity coefficients $\mu$.

\subsection{Viscous shock profile} \label{ssec.ViscShock}
\begin{figure}[!htbp]
	\begin{center}
		\begin{tabular}{cc} 
			\multicolumn{2}{c}{\includegraphics[width=0.9\textwidth]{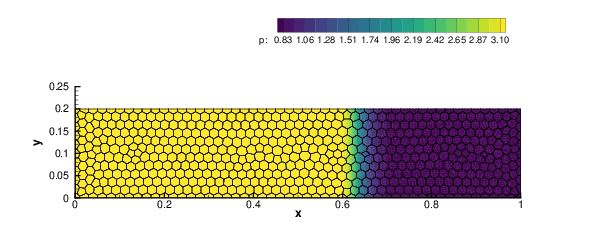}} \\
			\includegraphics[width=0.47\textwidth]{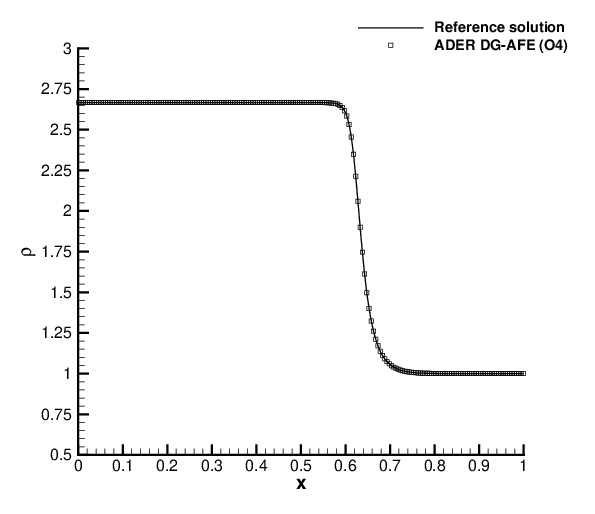} & 
			\includegraphics[width=0.47\textwidth]{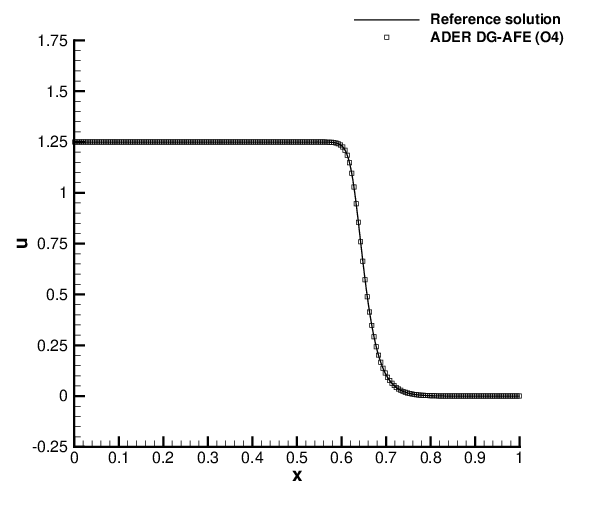} \\
			\includegraphics[width=0.47\textwidth]{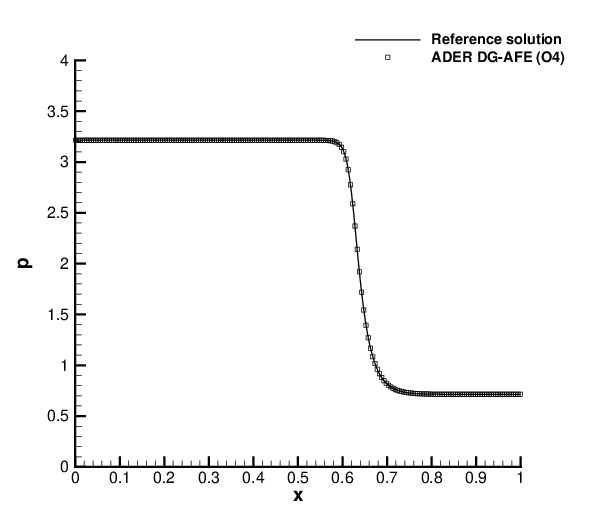} & 
			\includegraphics[width=0.47\textwidth]{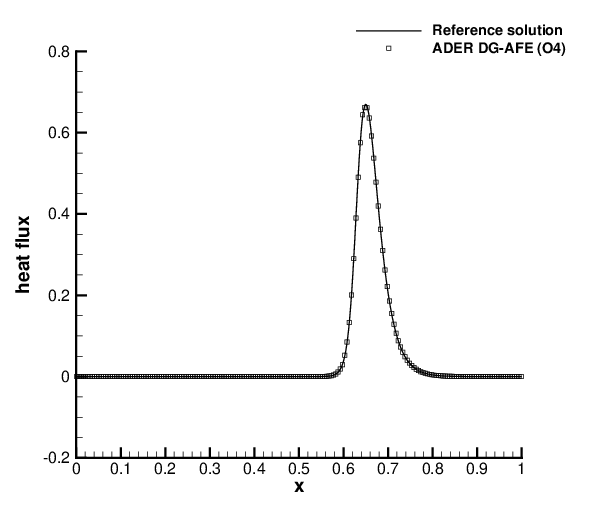} \\
		\end{tabular} 
		\caption{Viscous shock profile with shock Mach number $M_s=2$ and Prandtl number $Pr=0.75$ at time $t_f=0.2$. Top panel: Voronoi tessellation and pressure contours. Fourth order numerical solution with ADER DG-AFE scheme compared against the reference solution for density, horizontal velocity, pressure and heat flux (from middle left to bottom right panel): in particular, we show a one-dimensional cut of 200 equidistant points along the $x-$direction at $y=0.1$.}
		\label{fig.ViscousShock}
	\end{center}
\end{figure}

In order to verify the numerical method in the context of smooth supersonic viscous flows, we propose to solve the problem of an isolated viscous shock wave which is traveling into a medium at rest with a shock Mach number of $M_s > 1$, see also \cite{ADERNSE,ALEDG,BRMVCD21}. The analytical solution can be found in~\cite{Becker1923}, where the compressible Navier-Stokes equations are solved for the special case of a stationary shock wave at Prandtl number $Pr= 0.75$ with constant viscosity.
The dimensionless velocity of this stationary shock wave is $\bar u = \frac{u}{M_s \, c_0}$, with the upstream adiabatic sound speed $c_0=\sqrt{\gamma p_0/\rho_0}$. The exact solution for $\bar u$ can be computed by solving the following equation, see~\cite{Becker1923,ADERNSE}:
\begin{equation} 
\label{eqn.alg.u} 
\frac{|\bar u - 1|}{|\bar u - \kappa^2|^{\kappa^2}} = \left| \frac{1-\kappa^2}{2} \right|^{(1-\kappa^2)} 
\exp{\left( \frac{3}{4} \textnormal{Re}_s \frac{M_s^2 - 1}{\gamma M_s^2} x \right)},
\end{equation}
with
\begin{equation}
\kappa^2 = \frac{1+ \frac{\gamma-1}{2}M_s^2}{\frac{\gamma+1}{2}M_s^2}.
\end{equation}
Equation~\eqref{eqn.alg.u} allows the dimensionless velocity $\bar u$ to be obtained as a function of $x$. 
The form of the viscous profile of the dimensionless pressure $\bar p = \frac{p-p_0}{\rho_0 c_0^2 M_s^2}$ is given by
the relation 
\begin{equation}
\label{eqn.alg.p} 
\bar p = 1 - \bar u +  \frac{1}{2 \gamma}
\frac{\gamma+1}{\gamma-1} \frac{(\bar u - 1 )}{\bar u} (\bar u - \kappa^2).  
\end{equation}
Finally, the profile of the dimensionless density $\bar \rho = \frac{\rho}{\rho_0}$ is derived from the 
integrated continuity equation: $\bar \rho \bar u = 1$. A constant velocity field $u = M_s c_0$ is then superimposed to the solution of the stationary shock wave found in the previous steps, so that an unsteady shock wave traveling into a medium at rest is obtained.
The computational domain is the rectangular box $\Omega=[0,1]\times[0,0.2]$ which is discretized by a coarse grid made of a total number of unstructured cells $N_E=1120$, depicted in the top panel of Figure~\ref{fig.ViscousShock}. On the left side of the domain ($x=0$) the constant inflow velocity is prescribed, whereas outflow boundary conditions are imposed at $x=1$. Periodic boundaries are prescribed elsewhere. The initial condition involves a shock wave centered at $x=0.25$ propagating at Mach $M_s=2$ from left to right with a Reynolds number $Re=100$, thus the viscosity coefficient is set to $\mu=2\times 10^{-2}$. The upstream shock state is defined by $\mathbf{P}_0(\x,0)=(1,0,0,1/\gamma)$, hence $c_0=1$. The final time of the simulation is $t_f=0.2$ with the shock front located at $x=0.65$. Figure~\ref{fig.ViscousShock} shows a comparison of the fourth order numerical solution against the analytical solution at the final time, where an excellent matching is achieved. We compare the exact solution and the numerical solution for density, horizontal velocity component, pressure and heat flux $q_x=\kappa \, \frac{\partial T}{\partial x}$, with $T=p/(R\rho)$ being the temperature. We underline that this test case allows all terms contained in the Navier-Stokes system to be properly checked, since advection, thermal conduction and viscous stresses are present. Furthermore, despite the underlying one-dimensional structure of the exact solution, this problem actually becomes multidimensional due to the use of unstructured Voronoi meshes, where in general the edges of the control volumes are not aligned with the main flow field in the $x-$direction. 

The same simulation is also run using the fourth order version of the ADER DG-M scheme. Since the analytical solution is known, the errors in $L_2$ and $L_{\infty}$ norms for the main variables can be computed for both ADER DG-AFE and ADER DG-M methods at the final time of the simulation. Table~\ref{tab.viscousShock} reports the error analysis, which demonstrates the higher resolution achieved by the novel ADER DG-AFE schemes that exhibit much lower errors ($\approx 60\%$ less than the corresponding ADER DG-M errors) for all variables. 

\begin{table}[!htbp]  
	\caption{Error analysis for the viscous shock profile using both ADER DG-AFE and ADER DG-M schemes with fourth order of accuracy in space and time. The errors are measured in $L_2$ and $L_{\infty}$ norms and refer to the variables $\rho$ (density), horizontal velocity $u$ and pressure $p$ at the final time $t_{f}=0.2$.}  
	\begin{center} 
		\begin{small}
			\renewcommand{\arraystretch}{1.1}
			\begin{tabular}{c|cc|cc|cc}
				\multirow{2}{*}{Scheme} & \multicolumn{2}{c|}{Density ($\rho$)} & \multicolumn{2}{c|}{Velocity ($u$)} & \multicolumn{2}{c}{Pressure ($p$)} \\
				& $L_2$ & $L_{\infty}$ & $L_2$ & $L_{\infty}$ & $L_2$ & $L_{\infty}$ \\ 
				\hline
				ADER DG-AFE & 3.389E-05 & 1.518E-03 & 4.600E-05 & 1.757E-03 & 1.287E-04 & 5.450E-03 \\
				ADER DG-M\phantom{FE} & 8.692E-05 & 4.085E-03 & 1.337E-04 & 5.882E-03 & 3.687E-04 & 1.676E-02 \\
			\end{tabular}
		\end{small}
	\end{center}
	\label{tab.viscousShock}
\end{table}

\subsection{2D Taylor-Green vortex} \label{ssec.TGV}
\begin{figure}[!htbp]
	\begin{center}
		\begin{tabular}{cc} 
			\includegraphics[width=0.47\textwidth]{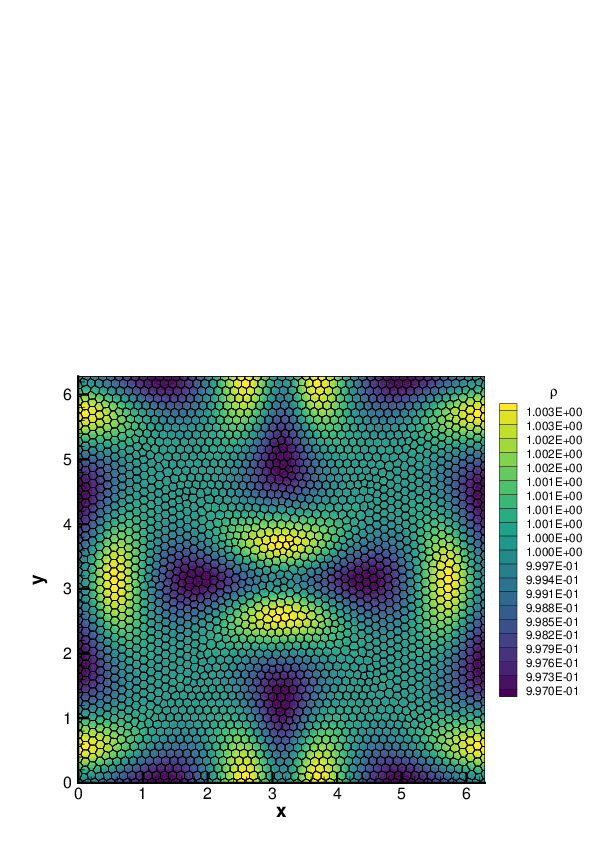} & 
			\includegraphics[width=0.47\textwidth]{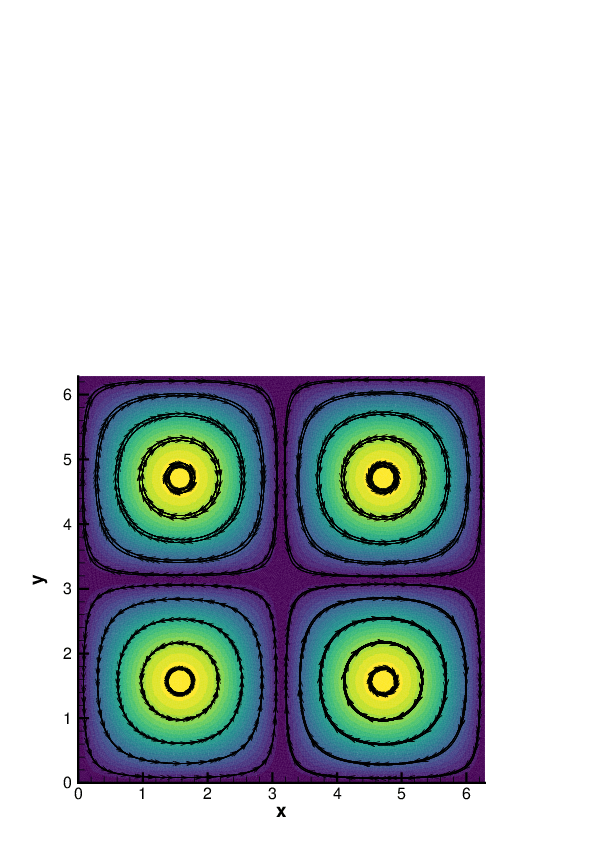} \\
			\includegraphics[width=0.47\textwidth]{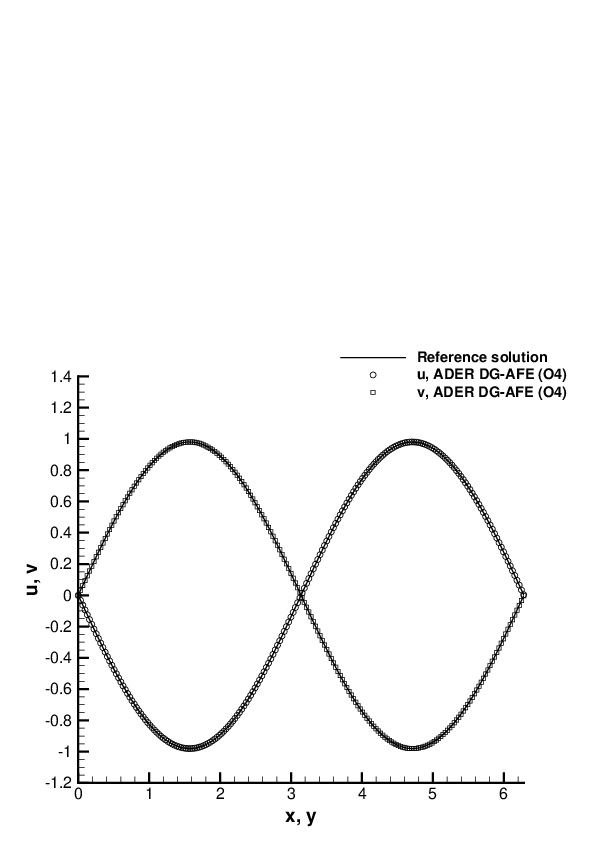} & 
			\includegraphics[width=0.47\textwidth]{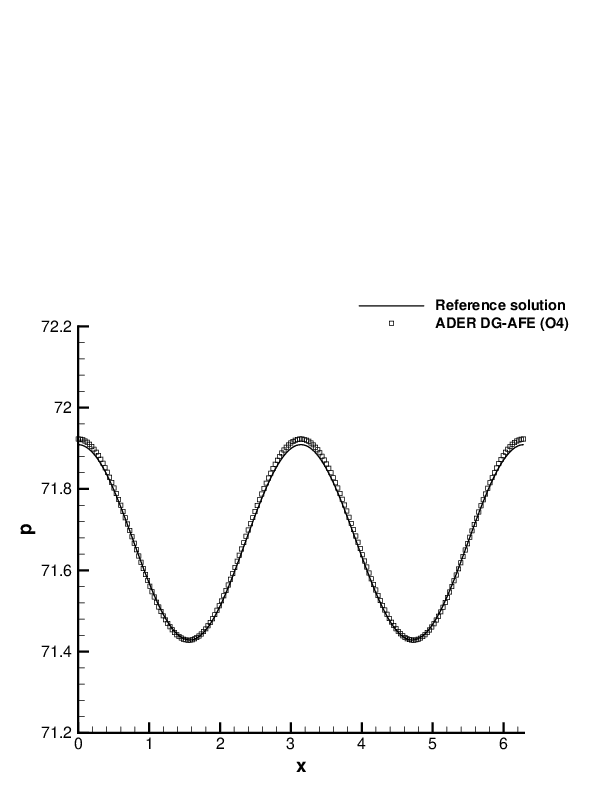} \\
		\end{tabular} 
		\caption{2D Taylor-Green vortex at time $t_f=1$ with viscosity $\mu=10^{-2}$. Exact solution of the Navier-Stokes equations and fourth order numerical solution with ADER DG-AFE scheme. Top: mesh configuration with density distribution (left) and vorticity magnitude with stream-traces (right). Bottom: one-dimensional cut of 200 equidistant points along the $x$-axis and the $y-$axis for the velocity components $u$ and $v$ (left) and for the pressure $p$ (right).}
		\label{fig.TGV2D}
	\end{center}
\end{figure}

A typical test problem used for the validation of numerical methods for the incompressible Navier-Stokes equations is the Taylor-Green vortex problem, for which an exact solution is known in two space dimensions:
\begin{eqnarray}
u(\x,t)&=&\phantom{-}\sin(x)\cos(y) \, e^{-2\nu t},  \nonumber \\
v(\x,t)&=&-\cos(x)\sin(y) \, e^{-2\nu t}, \nonumber \\
p(\x,t)&=& C + \frac{1}{4}(\cos(2x)+\cos(2y)) \, e^{-4\nu t},
\label{eq:TG_ini}
\end{eqnarray}
where the kinematic viscosity is given by $\nu=\frac{\mu}{\rho}$ with $\mu=10^{-2}$. To approach a low Mach regime, the additive constant for the pressure field is given by $C=100/\gamma$ and the density is initially set to $\rho(\x,0)=1$, while heat conduction is neglected, thus setting $\kappa=0$. The computational domain is defined by $\Omega(0)=[0;2\pi]^2$ with periodic boundaries imposed on each side and it is discretized with a total number of $N_E=2916$ Voronoi cells. The fourth order accurate numerical results are depicted in Figure~\ref{fig.TGV2D} at the final time of the simulation $t_f=1.0$. A very good agreement between the ADER DG-AFE scheme in the low Mach number regime and the exact solution of the incompressible Navier-Stokes equations can be observed, both for velocity and pressure.  We also plot the stream-traces of the velocity field and of the density perturbations.

\subsection{Compressible mixing layer} \label{ssec.MixLayer}
\begin{figure}[!htbp]
	\begin{center}
		\begin{tabular}{cc} 
			\includegraphics[trim= 5 5 5 5, clip, width=0.9\textwidth]{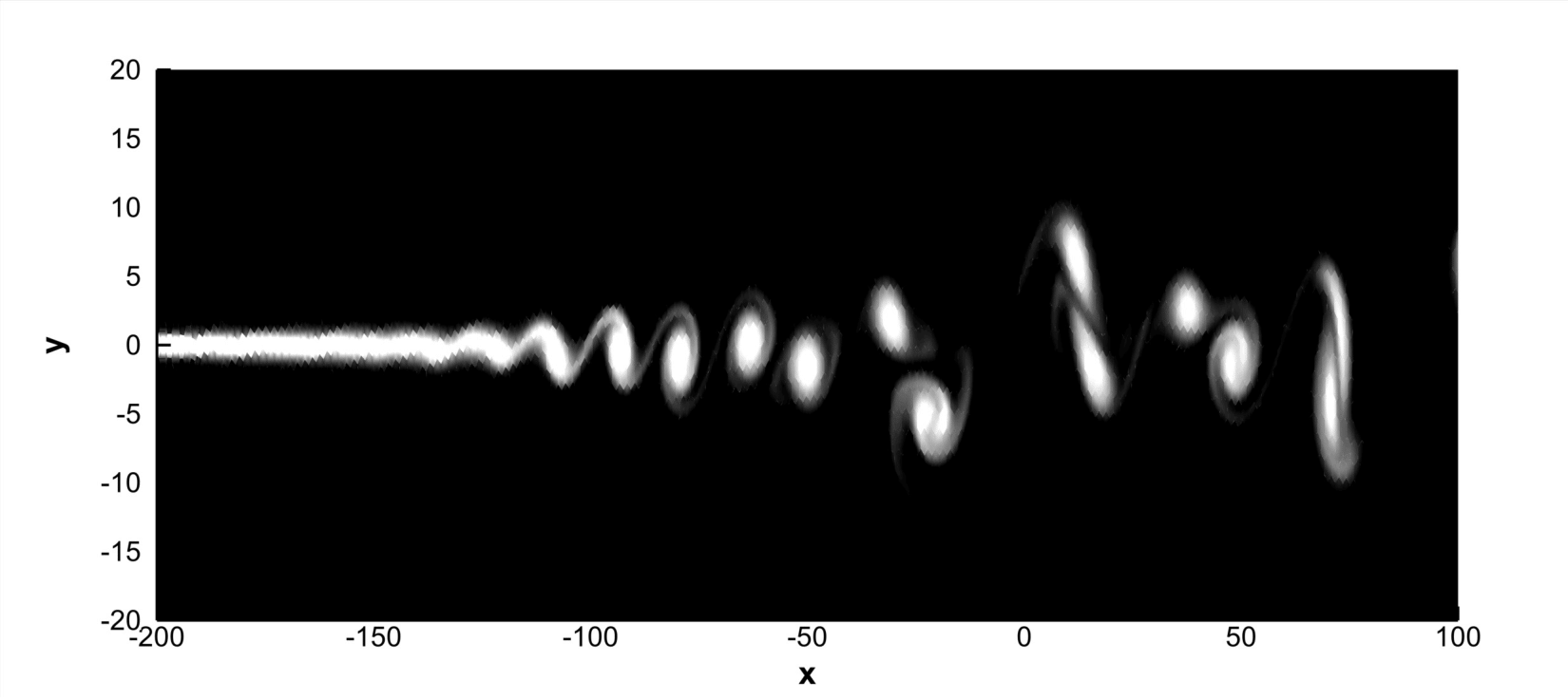}\\				
			\includegraphics[trim= 5 5 5 5, clip, width=0.9\textwidth]{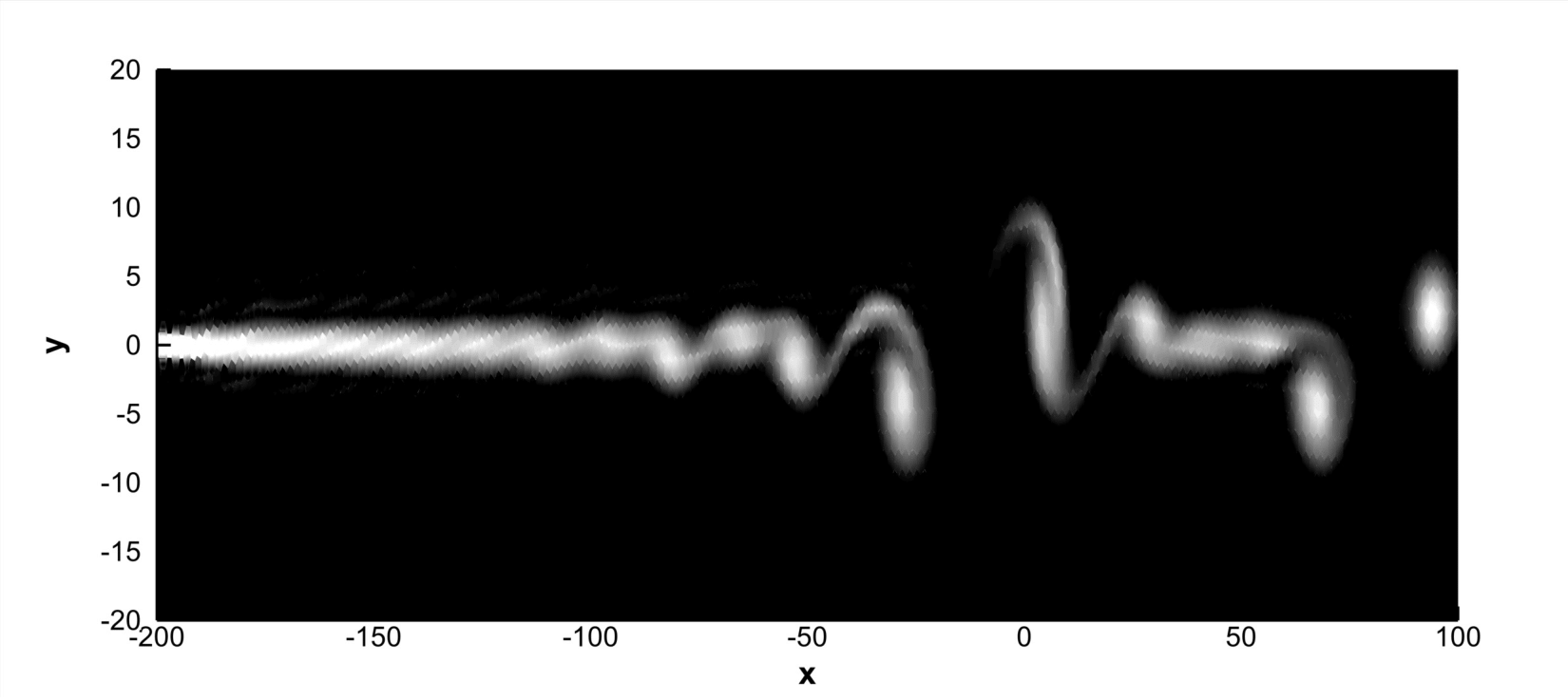}\\				
		\end{tabular} 
		\caption{Compressible mixing layer at time $t=500$. Third order numerical results for vorticity magnitude with ADER DG-AFE scheme (top row) and ADER DG-M scheme (bottom row). 40 contour levels in the range $[0.01;0.1]$ have been used for plotting the vorticity distribution.}
		\label{fig.MixingLayer500}
	\end{center}
\end{figure}

\begin{figure}[!htbp]
	\begin{center}
		\begin{tabular}{cc} 
			\includegraphics[trim= 5 5 5 5, clip, width=0.9\textwidth]{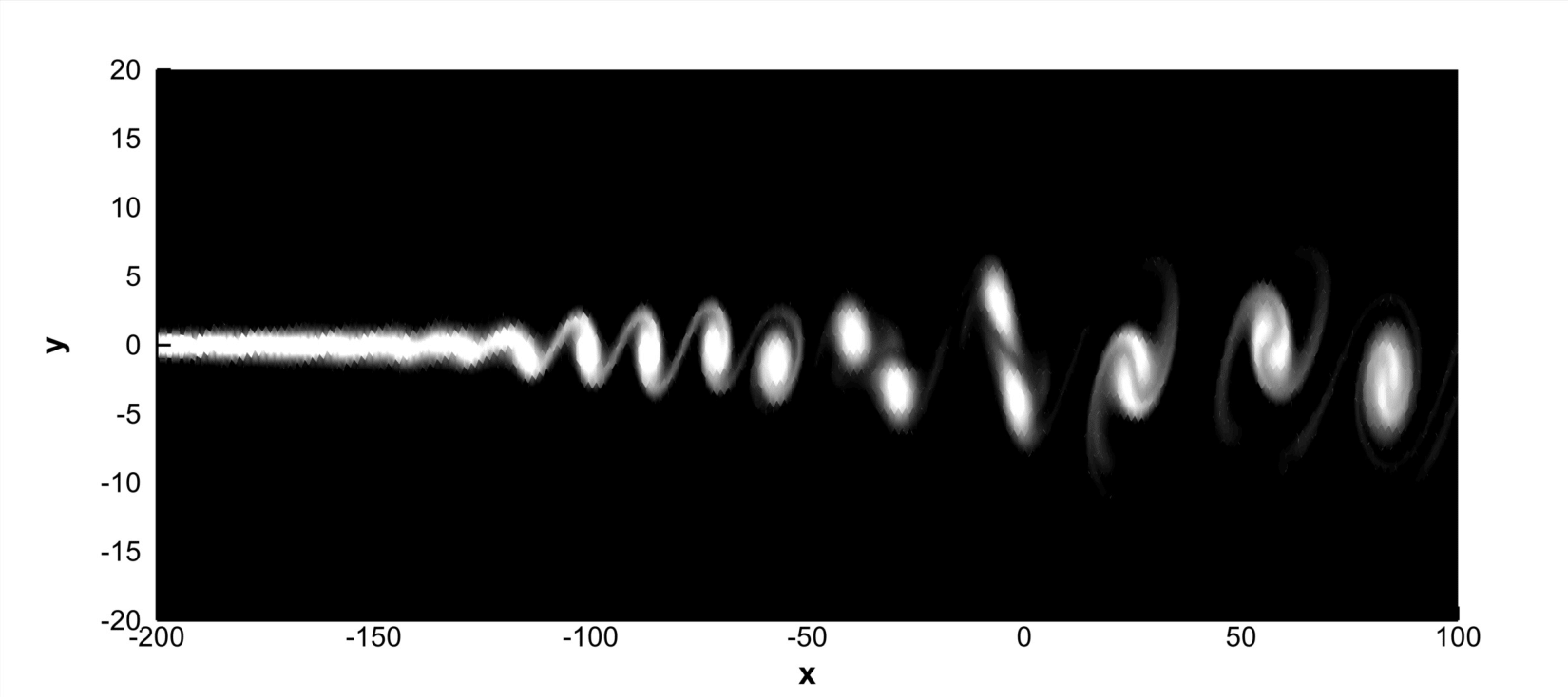}\\
			\includegraphics[trim= 5 5 5 5, clip, width=0.9\textwidth]{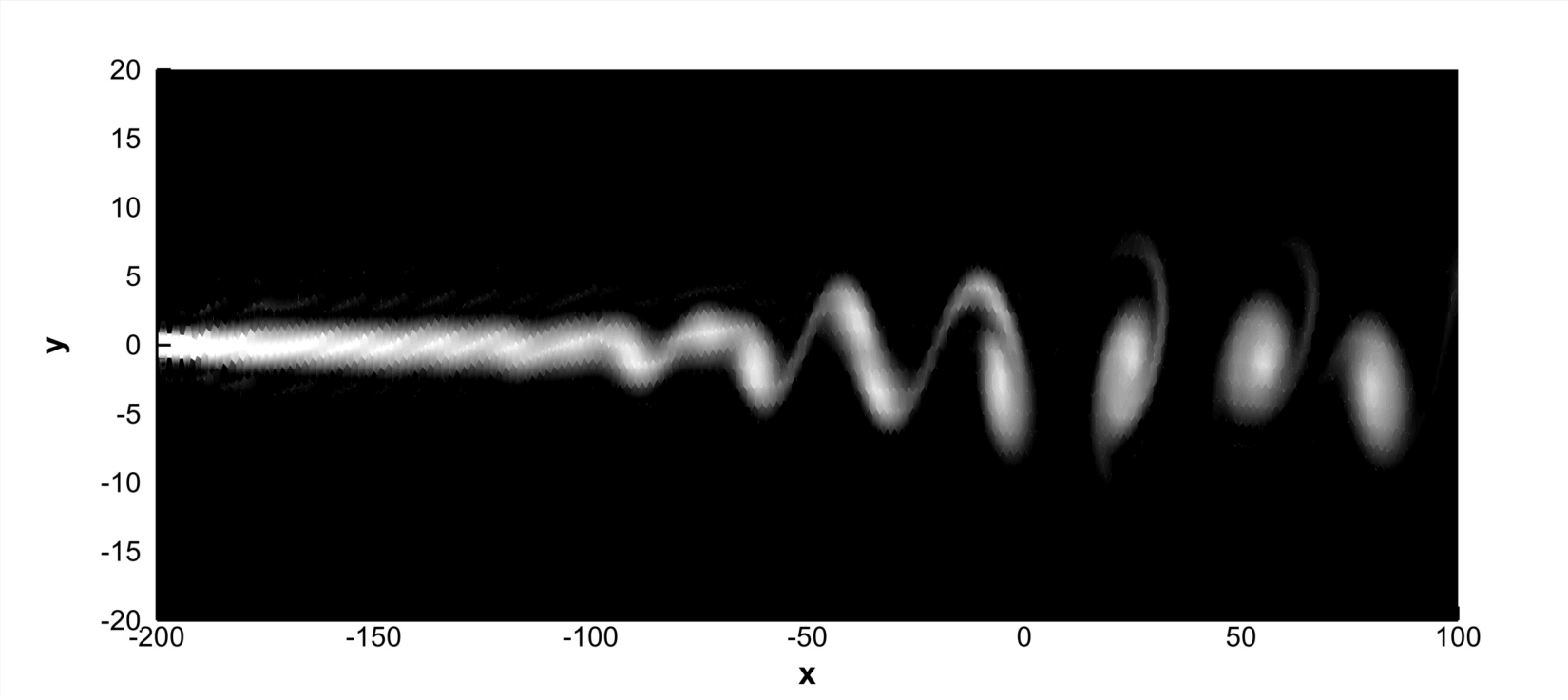}\\
		\end{tabular} 
		\caption{Compressible mixing layer at the final time $t=1596.8$. Third order numerical results for vorticity magnitude with ADER DG-AFE scheme (top row) and ADER DG-M scheme (bottom row). 40 contour levels in the range $[0.01;0.1]$ have been used for plotting the vorticity distribution.}
		\label{fig.MixingLayerEnd}
	\end{center}
\end{figure}

As last test to be presented in this work we consider an unsteady compressible mixing layer, originally studied by Colonius et al. in~\cite{Colonius}. The computational domain is given by $\Omega=[-200;200] \times [-50;50]$ which is paved with a total number of $N_E=15723$ unstructured cells. This problem involves two fluid layers moving with different velocities along the $x-$direction, thus for $y \to + \infty$ and $y \to -\infty$ we impose the free stream velocities $u_{+\infty}=0.5$ and $u_{-\infty}=0.25$, respectively. A smooth transition between the two velocities is initially imposed, thus the fluid state at $t=0$ is given according to~\cite{GPR1} as
\begin{equation}
\rho(\x,0) = \rho_0= 1, \quad \v(\x,0) = \v_0=(u_0,v_0)=\left( \begin{array}{c}
\frac{1}{8} \tanh(2y) + \frac{3}{8} \\ 0
\end{array} \right), \quad p(\x,0) = p_0 = \frac{1}{\gamma}.
\end{equation} 
The vorticity thickness $\chi$ at the inflow, with respect to which all lengths are made dimensionless, and the associated Reynolds number $Re_{\chi}$ write
\begin{equation}
\chi = \frac{u_{+\infty}-u_{-\infty}}{\max \left(\left.\frac{\partial u}{\partial y} \right|_{x=0}\right)}:=1, \qquad  Re_{\chi}=\frac{\rho_0 \, u_{+\infty} \, \chi}{\mu} = 500,
\end{equation}
with the viscosity coefficient set to $\mu=10^{-3}$, while heat	conduction is again neglected ($\kappa=0$). As done in~\cite{GPR1}, a perturbation $\delta(y,t)$ is introduced at the inflow boundary for all the variables, hence prescribing the following boundary condition at $x=0$:
\begin{equation}
\rho(0,y,t) = \rho_0 + 0.05 \, \delta(y,t), \quad \v(0,y,t) = \v_0 + \left( \begin{array}{c}
1.0 \\ 0.6
\end{array} \right)\, \delta(y,t), \quad p(0,y,t) = p_0 + 0.2 \, \delta(y,t).
\end{equation}
The periodic function $\delta(y,t)$ is given by
\begin{equation}
\delta(y,t) = -10^{-3} \exp(-0.25 y^2) \left[ \cos(\omega t) + \cos\left(\frac{1}{2}\omega t -0.028\right) + \cos\left(\frac{1}{4}\omega t +0.141\right) + \cos\left(\frac{1}{8}\omega t +0.391\right) \right],
\end{equation}
with the fundamental frequency of the mixing layer $\omega = 0.3147876$. Transmissive boundaries are imposed elsewhere. The final time is $t_f=1596.8$ and the simulation is carried out using the third order version of both ADER DG-AFE and ADER DG-M schemes. Figures~\ref{fig.MixingLayer500}-\ref{fig.MixingLayerEnd} show the vorticity magnitude of the flow field at $t=500$ and $t=t_f$ for both simulations, highlighting the differences between the two schemes. The novel ADER DG-AFE method using the continuous finite element subgrid basis clearly provides much better resolved flow patterns compared to the ADER DG-M method relying on the simple modal basis, and our new method is at the same time almost two times faster according to the profiling analysis reported in Table~\ref{tab.efficiency}.  

\section{Conclusions} \label{sec.concl}
In this article we have presented a new discontinuous Galerkin (DG) finite element method for the solution of nonlinear systems of conservation laws on unstructured polygonal Voronoi meshes in which the numerical solution is represented by a novel \textit{nodal basis}, making use of \textit{continuous finite element} basis functions defined on a subgrid \textit{within} each Voronoi element. 
The key novelty introduced in this paper is the use of \textit{piecewise polynomials} of degree $N$ as basis and test functions within each high order polygonal DG element, rather than using simply \textit{polynomials} of degree $N$ inside each element, as for example in the case of classical modal basis functions or for nodal basis functions based on the minimum necessary number of nodes \cite{GassnerPoly}. As usual in the context of DG schemes, also in our new approach the discrete solution is allowed to jump across element boundaries. The subgrid used within each polygon in order to define the subgrid finite element basis, also denoted by Agglomerated Finite Element (AFE) basis in this paper, has been defined by the sub-triangles connecting the barycenter with the vertices of each polygonal cell. 
Each triangular subcell can then be easily mapped to a universal reference triangle, hence allowing \textit{universal} mass, flux and stiffness matrices to be computed in the reference space. The basis functions are continuously extended by zero outside the subcell and only the inverse of the time stiffness matrix and the element mass matrix must be stored for each polygonal cell of the computational grid, whereas the remaining integrations can be \textit{efficiently} carried out by assembling the contributions of each subcell in a finite element fashion.   
Next, the new algorithm achieves also high order of accuracy in time thanks to the design of a space-time predictor-corrector approach relying on the ADER methodology.  
In particular, at the end of the predictor stage, the space-time polynomials representing the predictor solution are integrated in time. Then, a fully discrete one-step corrector DG scheme is applied, completing the integration in space by means of a variational formulation of the of the governing PDE.  

Thus, we have devised a novel \textit{quadrature-free} nodal DG scheme that is able to achieve high order of accuracy in space and time on general unstructured polygonal meshes. It has been carefully validated in terms of accuracy, robustness and computational efficiency by performing a suite of academic benchmarks for the compressible Euler and Navier-Stokes equations. The results have been compared against available exact or numerical reference solutions and both the performance and the accuracy of the new ADER DG-AFE schemes have been compared with the corresponding ADER DG methods based on classical modal basis functions defined in the physical space, showing the superior resolution capabilities and computational efficiency of the new algorithm proposed in this paper.

In future research we plan to extend the present scheme to three space dimensions, hence involving general polyhedral computational meshes with a subgrid made of tetrahedra. Within the context of all Mach number flows, a semi-implicit version of the DG scheme based on an explicit discretization of the nonlinear convective terms and an implicit treatment of the pressure and of the viscous terms will also be considered in the future, following the ideas presented in~\cite{TavelliCNS,BTBD20,BDT_cns}. This should also mitigate the time step restriction induced by the artificial viscosity limiter. More semi-implicit DG schemes for compressible flows and shallow water flows that may take advantage of the new nodal basis proposed in this paper can be found, e.g., in  \cite{Dolejsi1,Dolejsi2,Dolejsi3,GR10,TBR13,TumoloBonaventura2015}. 
Last but not least, motivated by the results obtained in~\cite{Klingenberg2015,KlingenbergPuppo,gaburro2018well,gaburro2018diffuse,Li_WBDG2018} for finite volume schemes applied to the Euler equations with gravity, and following the seminal approach of \cite{Castro2008} we also plan to introduce some well balanced structure preserving techniques in the new discontinuous Galerkin schemes presented in this paper, exploiting in particular the improved computational efficiency and the additional resolution capability of the new AFE basis functions in order to design new high order accurate and exactly well balanced DG schemes on polygonal meshes.

\section*{Acknowledgments}
W.~B. and M.~D. would like to thank the Italian Ministry of Instruction, University and Research (MIUR) to support this research with funds coming from PRIN Project 2017 (No. 2017KKJP4X entitled \quotew{Innovative numerical methods for evolutionary partial differential equations and applications}). W.~B. and M.~D. are both members of the INdAM GNCS group and E.~G is member of the CARDAMOM group at the Inria center of the University of Bordeaux. 
E.~G. gratefully acknowledges the support received from the European Union’s Horizon 2020 Research and Innovation Programme under the Marie Skłodowska-Curie Individual Fellowship \textit{SuPerMan}, grant agreement No. 101025563. 
The authors would also like to acknowledge the nice atmosphere and overall support provided by the Shark-FV 2019 workshop, \url{https://shark-fv.eu/old-conferences/shark2019}, where the main ideas contained in this paper were conceived and implemented.    



\bibliographystyle{plain}
\bibliography{biblioB}	

\end{document}